\documentclass[a4paper,12pt]{article}
\usepackage{amsmath,amssymb,amsfonts,graphicx}
\usepackage{subcaption}
\usepackage{epsfig}
\usepackage{mathtools}
\usepackage{color}
\usepackage{wrapfig}
\usepackage{comment}
\usepackage{ulem}
\usepackage{mathdots}
\usepackage{graphics}
\usepackage{epstopdf}
\usepackage{mathdots}
\usepackage{epsfig}
\usepackage{verbatim}
\usepackage{color}
\usepackage{array}
\usepackage{url}

\usepackage{mathrsfs}
\usepackage{eufrak}

\def\IR{{\mathbb R}}
\def\IC{{\mathbb C}}

\def\IL{{\mathbb L}}

\newcommand{\sIL}{{{{\mathbb L}_s}}}

\newcommand{\bA}{{\bf A}}
\newcommand{\bB}{{\bf B}}
\newcommand{\bC}{{\bf C}}
\newcommand{\bD}{{\bf D}}
\newcommand{\bE}{{\bf E}}

\newcommand{\bG}{{\bf G}}

\newcommand{\bY}{{\bf Y}}

\newcommand{\bL}{{\bf L}}
\newcommand{\bM}{{\bf M}}

\newcommand{\bH}{{\bf H}}

\newcommand{\bW}{{\bf W}}
\newcommand{\bR}{{\bf R}}

\newcommand{\bX}{{\bf X}}

\newcommand{\bV}{{\bf V}}

\newcommand{\bfa}{{\bf a}}
\newcommand{\bc}{{\bf c}}

\newcommand{\bs}{{\bf s}}
\newcommand{\bn}{{\bf n}}
\newcommand{\bd}{{\bf d}}

\newcommand{\br}{{\bf r}}

\newcommand{\bh}{{\bf h}}

\newcommand{\bv}{{\bf v}}
\newcommand{\bw}{{\bf w}}
\newcommand{\bz}{{\bf z}}

\newcommand{\bee}{{\bf e}}
\newcommand{\bff}{{\bf f}}

\newcommand{\cO}{ {\cal O} }

\newcommand{\bphi}{ \boldsymbol{\phi} }

\newcommand{\bomega}{ \boldsymbol{\omega} }
\newcommand{\bOmega}{ \boldsymbol{\Omega} }

\newcommand{\bell}{\boldsymbol{\ell}}

\newcommand{\bLambda}{\boldsymbol{\Lambda}}

\newcommand{\cR}{ {\cal R} }

\title{Case study: Approximations of the Bessel Function\\[.5mm]}
\author{D.S. Karachalios\footnotemark[1]~ , I.V. Gosea \footnotemark[1]~, Q. Zhang\footnotemark[2]~, A.C. Antoulas \footnotemark[2]~.
}

\begin{document}
\maketitle
\renewcommand{\thefootnote}{\fnsymbol{footnote}}

\footnotetext[1]{
Data-Driven System Reduction and Identification Group, Max Planck Institute for Dynamics of Complex Technical Systems, Magdeburg, Germany e-mail:
{\tt karachalios , gosea@mpi-magdeburg.mpg.de}
}
\footnotetext[2]{
   Department of Electrical and Computer Engineering, 
   Rice University,Houston, e-mail:  
   {\tt qz18, aca@rice.edu}
   }
\begin{abstract}
The purpose of this note is to compare various approximation methods as applied to the inverse of the Bessel function of the first kind $\frac{1}{J_{0}(s)}$, in a given domain of the complex plane.
\end{abstract}
\section{Introduction}
The Bessel function of the first kind and order $n\in\mathbb{N}$, is defined as
\begin{equation}
J_{n}(s)=\frac{1}{2\pi{i}}\oint e^{(\frac{s}{2})(t-\frac{1}{t})} t^{-n-1}dt.
\end{equation}
In the sequel we will consider only the case  $n=0$. Our goal is to approximate 
\begin{equation}
H(s)=\frac{1}{J_{0}(s)}, s\in \IC.
\end{equation}
In particular we will be interested in the approximation of $H(s)$ in the rectangle 
\begin{equation}
\bOmega=[0,10]\times[-1,1]\subset\IC,
\end{equation}
of the complex plane. The magnitude of $H(s)$, in $\bOmega$, is as follows:
\begin{figure}[h]  
\begin{center}
\includegraphics[width=0.6\linewidth, height=3cm]{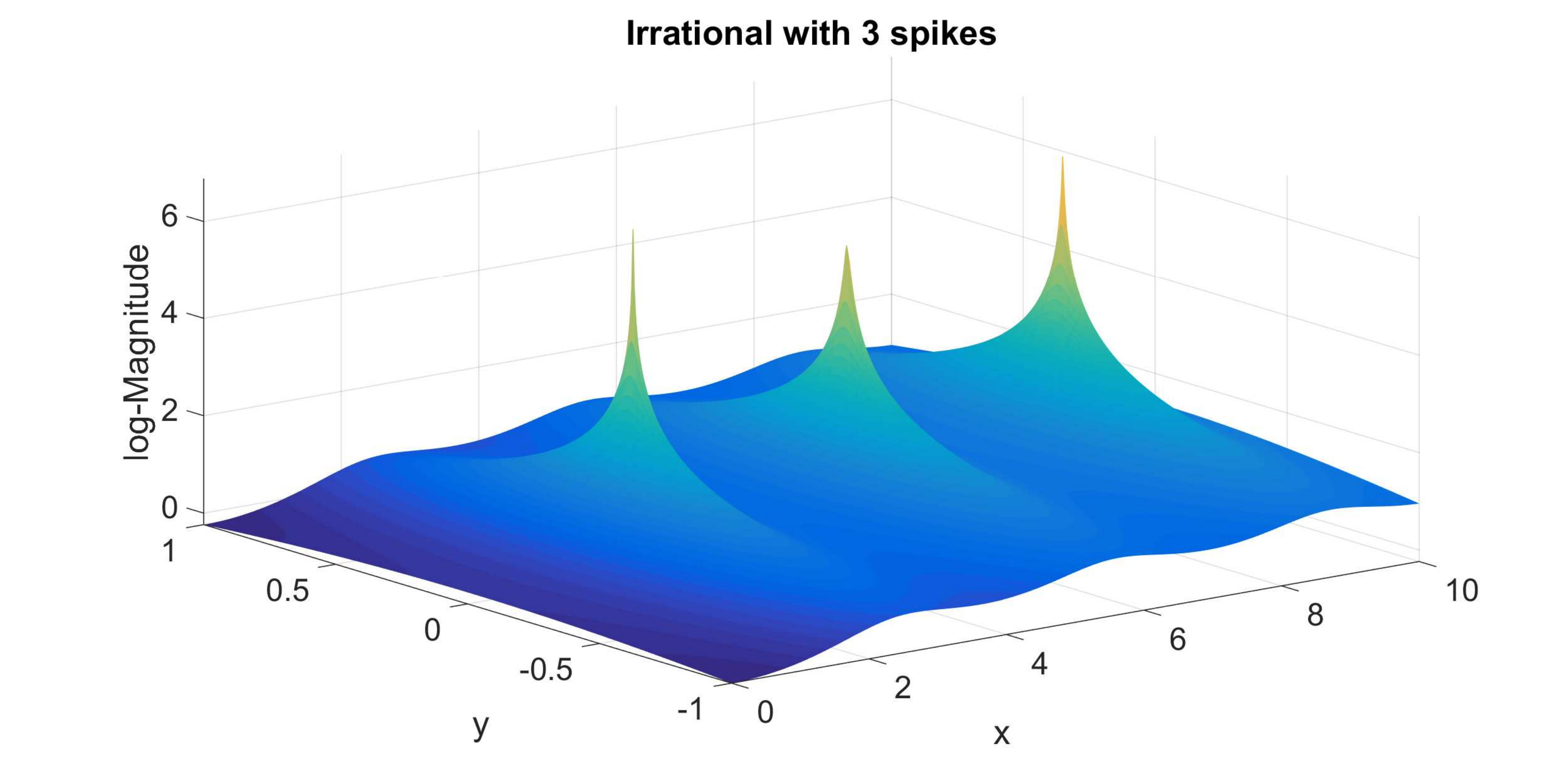}
\caption{$ln|H(s)|$, with $s\in\bOmega=[0,10]\times[-1,1]$}.
\end{center} 
\end{figure}
\newpage

The following interpolatory reduction methods will be used:
\begin{enumerate}
\item The Loewner Framework and (recursive) \cite{tutorial},
\item the AAA approach \cite{AAA},
\item Vector Fitting \cite{VF,DrmacVF}.
\end{enumerate}
\subsection{The Loewner Framework} \label{sec:Loewner}
The Loewner framework is a data-driven reduction method. It uses measured or computed data produced by an underlying dynamical system. The construction of reduced linear models is based on rank-revealing factorization. Here is a brief description of this method. For details, see \cite{tutorial}. Given is a set of pairs of complex numbers:
\begin{center}
$\{(s_{k}, \bphi_{k}): k=1,...,N\}$,
\end{center}
with $s_{k}\in\IC$, $\bphi_{k} \in\IC^{p\times m}$. We partition the data in two disjoint sets:
\begin{itemize}
\item left data: $(\mu_{j},\bf {\bell_{j}}, \bv_j)$, $j=1,...,q$
\item right data: $(\lambda_{i},\br_{i}, \bw_i)$, $i=1,...,k$
\end{itemize}

The objective is to find $\bH(s) \in \IC^{p\times m}$ such that:
\begin{equation}
\bH(\lambda_{i})\br_{i}\approx\bw_{i}, \quad \bf {\bell_j^*}\bH(\mu_{j})\approx\bv_j^{*}.
\end{equation}
In this report, we consider approximation of a SISO (Single Input Single Output) system, i.e. $m=p=1$, in which case the directions $\bf {\bell_{j}}$ and $\br_{i}$ can be chosen equal to 1. The associated Loewner and shifted Loewner matrices $\IL$, $\IL_{s}$, are defined as:
\begin{equation}
\IL=\left[\begin{array}{ccc}
\frac{\bv_1-\bw_1}{\mu_1-\lambda_1} & \cdots &
\frac{\bv_1-\bw_k}{\mu_1-\lambda_k} \\
\vdots & \ddots & \vdots \\
\frac{\bv_q-\bw_1}{\mu_q-\lambda_1} & \cdots &
\frac{\bv_q-\bw_m}{\mu_q-\lambda_k} \\
\end{array}\right],
\sIL=\left[\begin{array}{ccc}
\frac{\mu_1\bv_1-\bw_1\lambda_1}{\mu_1-\lambda_1} & \cdots &
\frac{\mu_1\bv_1-\bw_k\lambda_k}{\mu_1-\lambda_k} \\
\vdots & \ddots & \vdots \\
\frac{\mu_q\bv_q-\bw_1\lambda_1}{\mu_p-\lambda_1} & \cdots &
\frac{\mu_q\bv_q-\bw_k\lambda_k}{\mu_p-\lambda_k} \\
\end{array}\right],
\end{equation}
while $(W)_{i}=w_{i}$, $(V)_{j}=v_{j}$, $i=1,...,q$, $j=1,...,k$. In the case of just enough (minimal amount of) data, $(\IL,\sIL)$, is regular. In this case $\{\hat{\bE}=-\IL,\hat{\bA}=-\sIL,\hat{\bB}=\bV,\hat{\bC}=\bW\}$, is a minimal realization of the data. Hence, $\bG(\bz)=\bW(\sIL-z\IL)^{-1}\bV$, is the required minimal solution of the problem. In the case of redundant data we use the rank revealing SVD decompositions
\begin{equation}
\begin{bmatrix} \IL , \sIL  \end{bmatrix}=\bY\Sigma_{\bell}\tilde{\bX}^*, \quad
\begin{bmatrix} \IL \\\sIL  \end{bmatrix}=\tilde{\bY}\Sigma_{\br}\bX^*.
\end{equation}
where $\Sigma_{\bell}$ and $\Sigma_{\br}$ are $r\times r$, with $r\ll N$. \\The quadruple $\{\hat{\bE}=-\bY^*\IL\bX,\hat{\bA}=-\bY^*\sIL\bX,\hat{\bB}=\bY^*\bV,\hat{\bC}=\bW\bX\}$, is then, the realization of an approximate data interpolant of order r.
\begin{equation}
\bH(s)=\hat{\bC}(s\hat{\bE}-\hat{\bA})^{-1}\hat{\bB}
\end{equation}  
\subsubsection{Computation of poles/zeros and projected interpolation points\\ Skelton \cite{Skelton}, Van Dooren \cite{MRMIMO}, Grimme \cite{GrimmeThesis}}
Having a rational representation in state-space form, $H(s)=\hat{\bC}(s\hat{\bE}-\hat{\bA})^{-1}\hat{\bB}$ is crucial for computing its poles (roots of the denominator) and zeros (roots of the numerator), since these key quantities can be obtained, for example, by computing the finite eigenvalues resulting from the following generalized eigenvalue problems that involve the state-space matrices \cite{Ionita}
\begin{itemize}
\item poles: $\sigma(\hat{\bA},\hat{\bE})$
\item zeros: $\sigma\Bigg(\Bigg[\begin{tabular}{ c c }
  $\hat{\bA}$ & $\hat{\bB}$  \\
  $\hat{\bC}$ & $\hat{\bD}$  \\ 
\end{tabular}\Bigg],\Bigg[\begin{tabular}{ c c }
  $\hat{\bE}$ & $\textbf{0}$  \\
  $\textbf{0}$ & $\textbf{0}$ \\ 
\end{tabular}\Bigg]\Bigg)$, where in this case $\hat\bD=\textbf{0}$.
\end{itemize}

\textbf{Interpolatory projection matrices - (Skelton/Grimme/Van Dooren)}\\
We introduce the following quantities which play the role of projectors in the sequel. For this we need the sets of the complex numbers namely $\mu_{i}, i=1,...,q$, and $\lambda_{j}, j=1,...,k$ which we will refer to as left and right interpolation points:
\begin{equation}
\cR=[\bphi(\lambda_1)\bB \dots \bphi(\lambda_k)\bB ]\in\IC^{n\times k},\quad \cO=[\bC\bphi(\mu_1) \dots \bC\bphi(\mu_q)]^T\in\IC^{q\times n}.
\end{equation}
These are called the generalized \textbf{controllability} and the generalized \textbf{observability} matrices. With $\bLambda=diag[\lambda_1,\dots,\lambda_k]$, $\bM=diag[\mu_1,\dots,\mu_q]$, $\textbf{e}_m=[1 \dots 1]^{T}\in\IR^{m}$. Then, $\cR$ and $\cO$ satisfy the Sylvester equations:
\begin{equation}
\bE\cR\bLambda-\bA\cR=\bB\textbf{e}_{k}^{T}, \quad \bM\cO\bE-\cO\bA=\textbf{e}_{q}\bC.
\end{equation}
By multiplying the above equations with $\cO,\cR$ we obtain:
\begin{center}
$\IL_s-\IL\bLambda=\bV\bR$ and $\IL_s-\bM\IL=\bL\bW$
\end{center}
By adding/subtracting appropriate multiples of these expressions it follows that the Loewner quadruple satisfies the Sylvester equations
\begin{equation}
\bM\IL-\IL\bLambda=\bV\bR-\bL\bW, \quad \bM\IL_s-\IL_s\bLambda=\bM\bV\bR-\bL\bW\bLambda.
\end{equation} 
Given the projectors $\bX,\bY\in\IC^{n\times k}$, let the reduced quantities be
\begin{equation}
\hat\IL=\bX^*\IL\bY, \quad\hat\IL_s=\bX^*\IL_s\bY,\quad \hat\bV=\bX^*\bV,\quad\hat\bL=\bX^*\bL,\quad\hat\bW=\bW\bY,\quad \hat\bR=\bR\bY.
\end{equation}
The associated $\hat\bLambda$, and $\hat\bM$, must satisfy the projected equations above, so,
\begin{equation}
\hat\IL_s-\hat\IL\hat\bLambda=\hat\bV\hat\bR,\quad \hat\IL_s-\hat\bM\hat\IL=\hat\bL\hat\bW.
\end{equation}
Let the associated generalized eigenvalue problem be:
\begin{equation}
\hat\bLambda=eig(\hat\IL_s-\hat\bV\hat\bR,\hat\IL),  \hat\bM=eig(\hat\IL_S-\hat\bL\hat\bW,\hat\IL).
\end{equation}
From the above equations we are able to compute the \textbf{projected interpolation} points.
\subsection{The recursive Loewner approach}
A brief description is as follows. For details see \cite{tutorial}\\
\textbf{First step:}\\
Randomly choose two measurements. Put one of them into left measurement set and put another one into right measurement set. Generate the resulting approximant.\\
\textbf{Second step:}
\begin{itemize}
\item Calculate the error between the approximant and the measurement values at the sample points.
\item Put the two measurements which has largest error into the left/right measurement sets respectively.
\item Construct the new approximant by means of the recursive Loewner framework. Continue until all measurements are used.
\end{itemize}
\subsection{The AAA algorithm} \label{sec:AAA}
This algorithm \cite{AAA} uses \textbf{barycentric representation} of interpolates :
\begin{equation}
\br(\bs)=\frac{\bn(\bs)}{\bd(\bs)}=\frac{\sum_{\j=1}^m \frac{\bomega_{j}\bff_{j}}{\bs-\bs_{j}}}{\sum_{\j=1}^m \frac{\bomega_{j}}{\bs-\bs_{j}}}.
\end{equation}

\textbf{First step:} \\
Choose randomly one measurement from the measurement set. Construct the first order system $m=1$, where $\bs_1$ is the sample point and $\bff_1$ is the measurement value.

\textbf{Second step:} \\
Define the error as $\bee(\bomega)=\sum_{j=2}^n(\bff_j\bd(\bs)-\bn(\bs))^2$, where $\bomega$, describes the weights. Using least square we can calculate $\bomega$.

\textbf{Third step:} \\
\textbf{for} $j=2,...,k$
\begin{itemize}
\item Calculate the error between approximant and measurement value at sample points. 
\item Interpolate the measurement which has the largest error to increase the order of barycentric representations $m=j$ such that $\bs_j$, is the sample point and $\bff_j$, is the measurement value. 
\item Using least squares, calculate $\bomega$.
\end{itemize}
\textbf{end}\\
Return the barycentric representation approximation: $\br(\bs)=\frac{\bn(\bs)}{\bd(\bs)}=\frac{\sum_{\j=1}^k \frac{\bomega_{j}\bff_{j}}{\bs-\bs_{j}}}{\sum_{\j=1}^k \frac{\bomega_{j}}{\bs-\bs_{j}}}$. For more details see \cite{AAA}.
\newpage
\subsection{Vector Fitting} \label{sec:VF}
Condider the rational function approximation
\begin{equation}
f(s)=\sum_{n=1}^r \frac{\bc_n}{\bs-\bfa_n}+\bd+\bs\bh.
\end{equation}
The residue $\bc_n$, and poles are real or complex conjugate pairs, while $\bd$ and $\bh$, are real. The problem is to estimate all coefficients in the above expression so that the least square error gets minimized. The above problem is a non linear problem, because $\bfa_n$, appears in the denominator. Vector fitting addresses this issue by solving a sequence of linear problems. The procedure is as follows. Fixing the quantities $\bar{\bfa_i}$, solve for $\bar{\bc}_i$, $\bc_i$, $\bd$ and $\bh$, by evaluating the expression
\begin{equation}
\sum_{n=1}^r\bigg(\frac{\bar{\bc}_{n}}{s-\bar{\bfa}_n}+1\bigg)f(\bs)=\sum_{n=1}^r\frac{\bc_n}{\bs-\bar{\bfa}_n}+\bd+\bs\bh.
\end{equation}
at each interpolation point and using least squares to minimize the error. The procedure is repeated until convergence. For details, see \cite{VF,DrmacVF} or \url{http://www.sintef.no/projectweb/vectfit/}.

\section{Approximation Example}
\textbf{Choosing the interpolation points}. The aim is to approximate the irrational function $H(s)$, in Figure 1, in the domain $\bOmega$, which is sampled in two different ways. First using a structured grid with 2121 interpolation points (left pane below), and second, using a uniformly distributed grid with 2000 points (right pane below).
\begin{figure}[h]
\begin{subfigure}{.5\textwidth}
  \centering
  \includegraphics[width=1\linewidth]{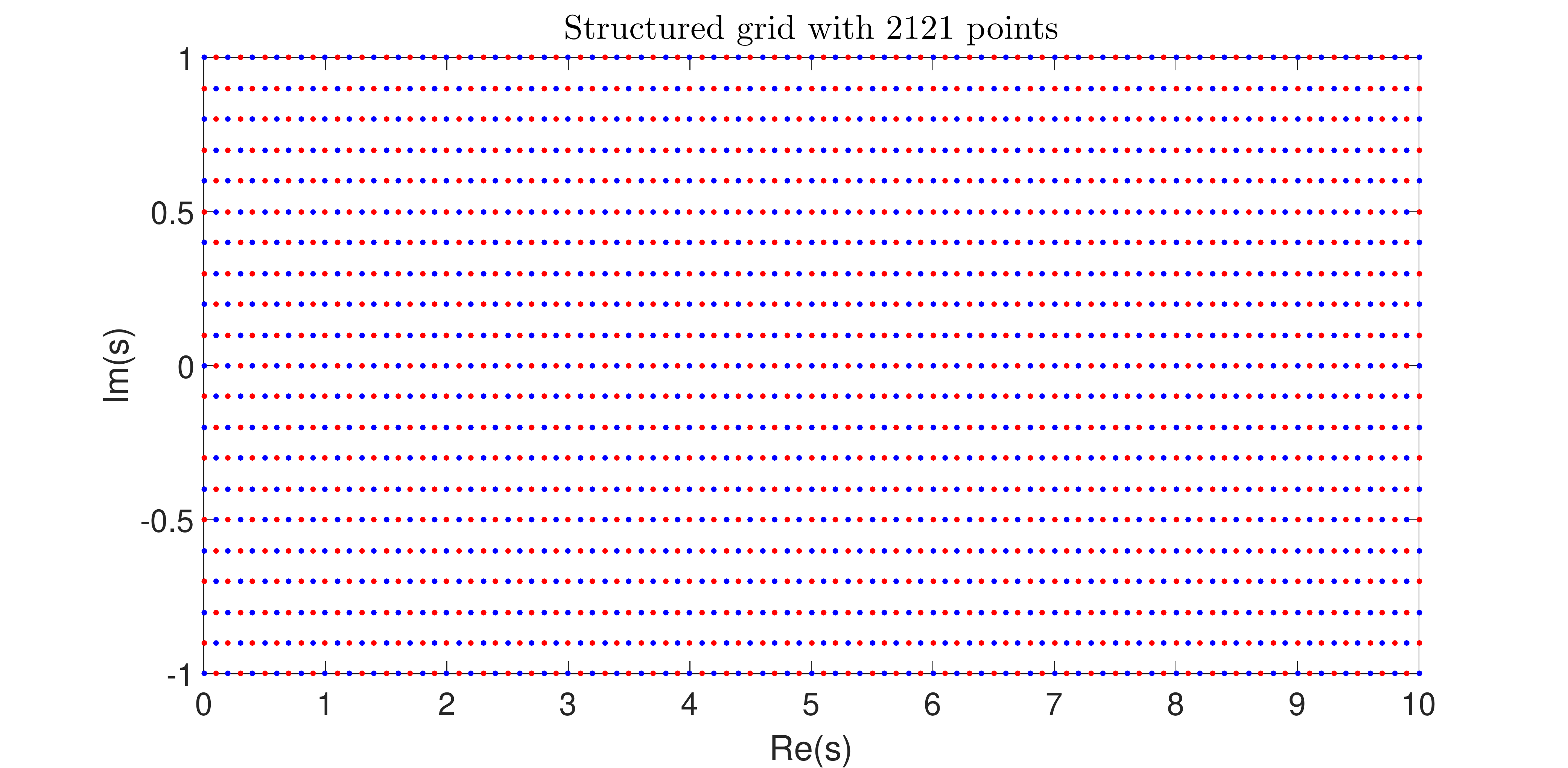}
  \caption{Structured grid with 2121 points}
  \label{fig:sfig1}
\end{subfigure}%
\begin{subfigure}{.5\textwidth}
  \centering
  \includegraphics[width=1\linewidth]{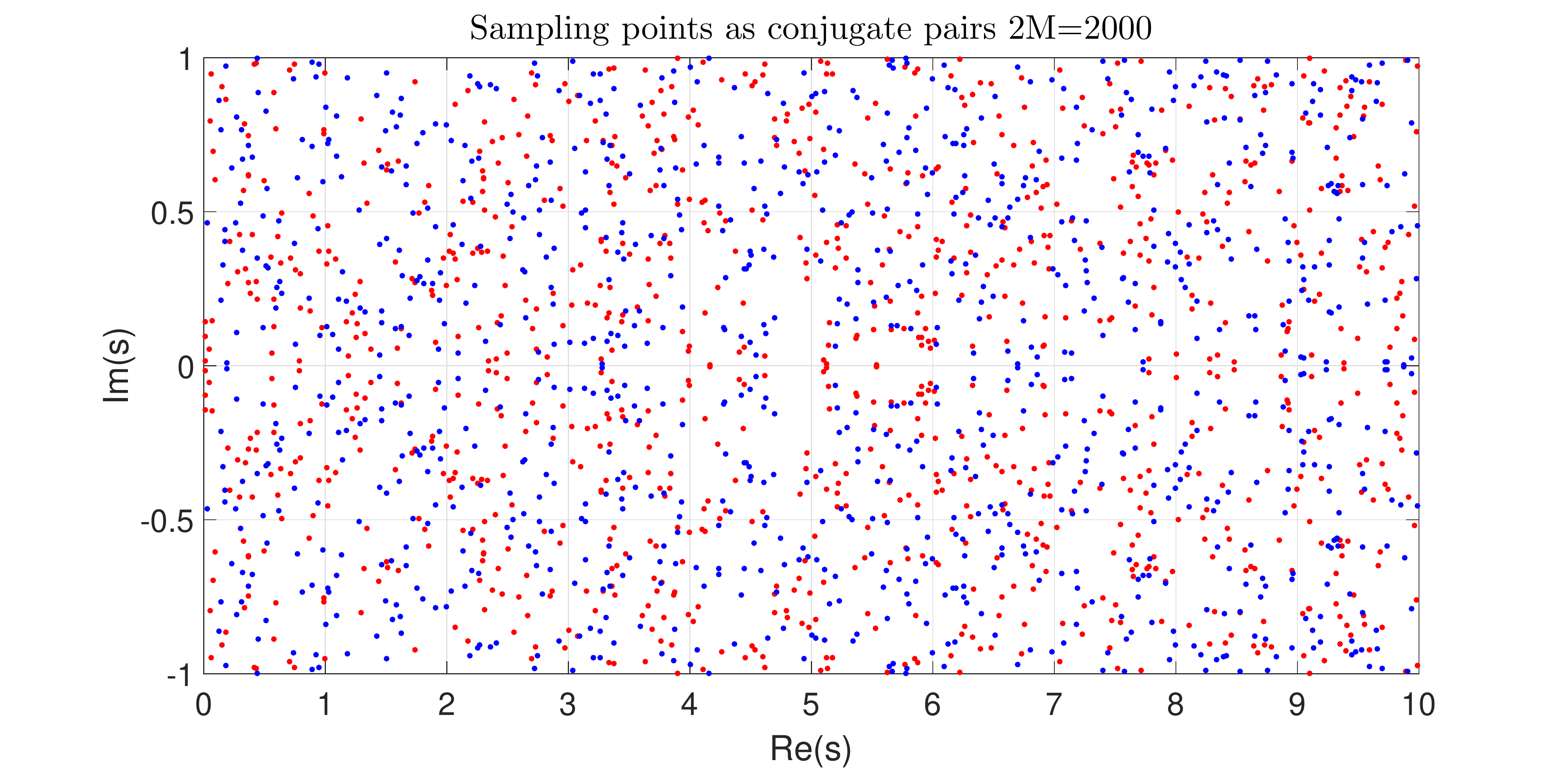}
  \caption{Uniformly grid with 2000 sampling points.}
  \label{fig:sfig2}
\end{subfigure}
\caption{In both cases we took under consideration conjugate pairs in order to enforce real symmetry.}
\label{fig:fig}
\end{figure} 
\subsection{Interpolation points as in Figure 2a}
\subsubsection{The Loewner Framework applied to $H(s)$}
First we will construct approximants $H_{r}(s)$, ($r$, order) of $H(s)$, using the interpolation points as shown in the left hand figure, section 2, page 5. Below is a description of these points with the corresponding samples of $H(s)$, (Figure 3a) and the (normalized) singular values of the augmented Loewner matrix (Figure 3b). We truncate at r=11 with corresponding accuracy $O(10^{-12})$.
\begin{figure}[h!]
\begin{subfigure}{.5\textwidth}
  \centering
  \includegraphics[width=1\linewidth]{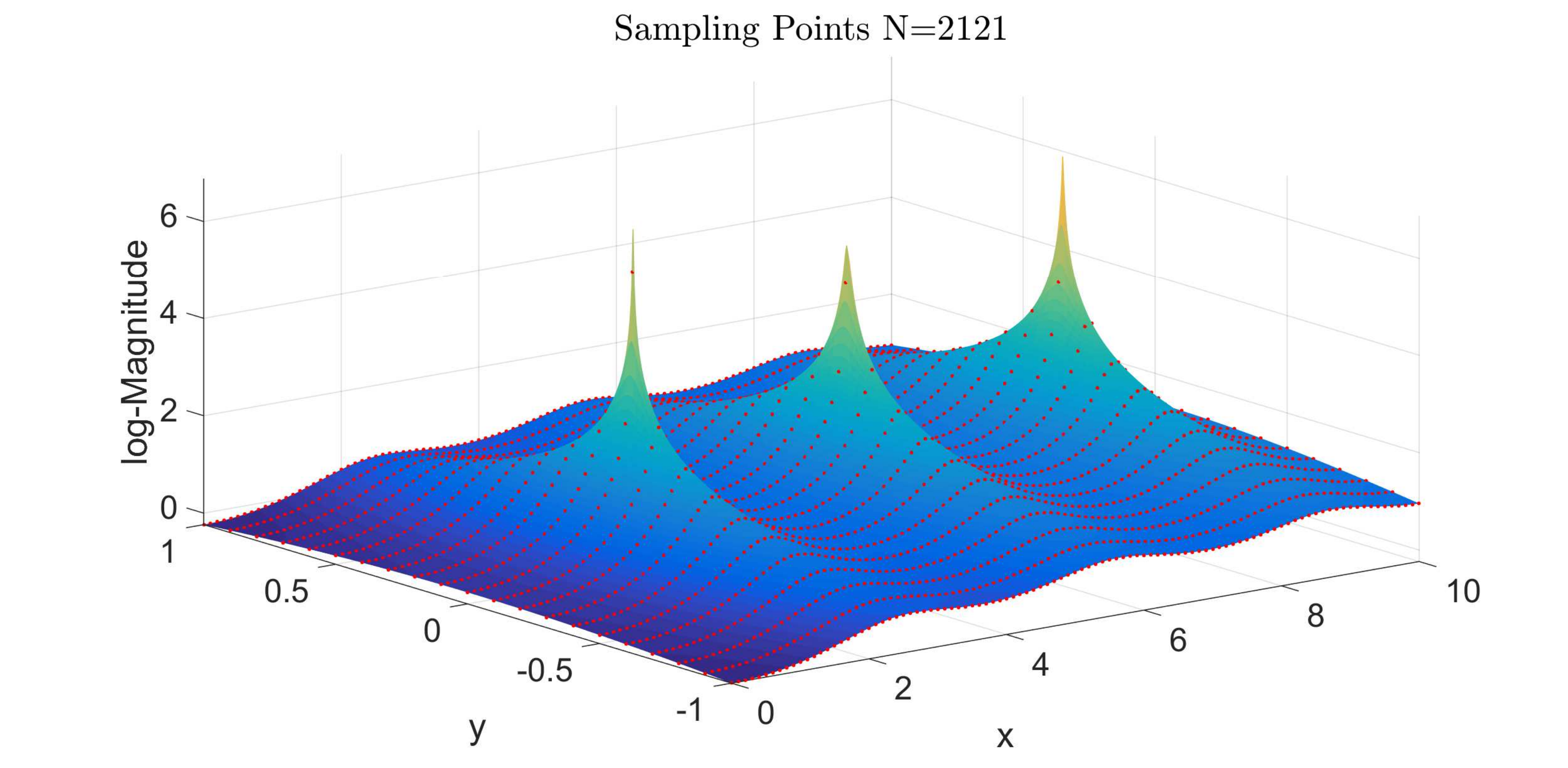}
  \caption{Sampling points over $\bOmega=[0,10]\times[-1,1]$},
  \label{fig:sfig1}
\end{subfigure}%
\begin{subfigure}{.5\textwidth}
  \centering
  \includegraphics[width=1\linewidth]{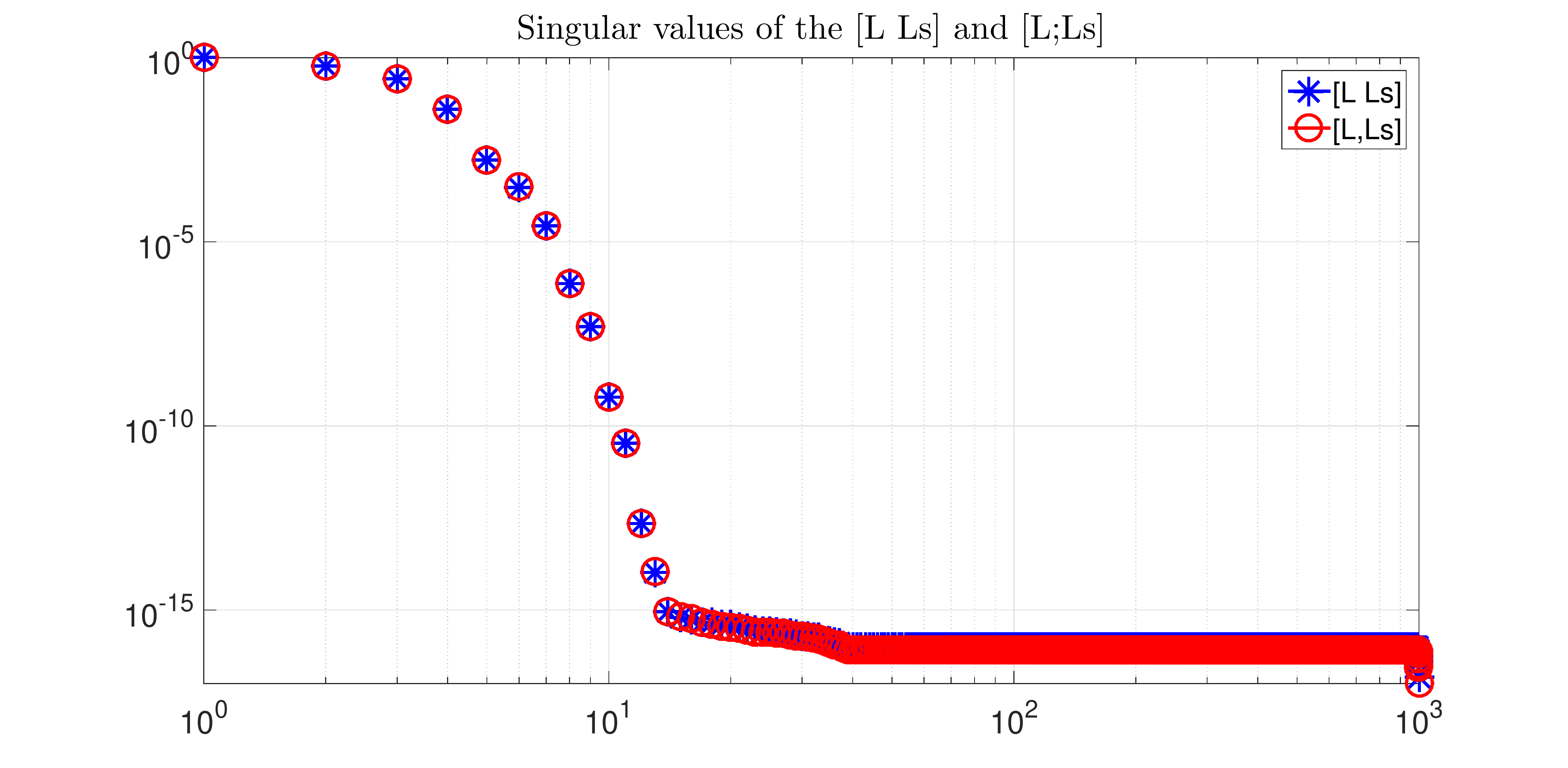}
  \caption{Fast decay of the singular values.}
  \label{fig:sfig2}
\end{subfigure}
\caption{Interpolation scheme of the irrational function and the corresponding rank reavealing SVD of the augmented Loewner matrix.}
\end{figure}
\newpage
The next plot (Figure 4) shows the 11th order approximant $H_r(s)$, superimposed on the plot of $H(s)$, in $\bOmega$, over the dense grid $\bOmega_{grid}=[x_{1},...,x_{500}]\times[y_{1},...,y_{500}]\subset \bOmega$. As a result we see the error plot in Figure 4b which shows for each point in the $\bOmega_{grid}$, the modulo $|H_{r}(s)-H(s)|$, with $s \in \bOmega_{grid}$. The order of error is $O(10^{-11})$.
  
\begin{figure}[h]
\begin{subfigure}{.5\textwidth}
  \centering
  \includegraphics[width=1\linewidth]{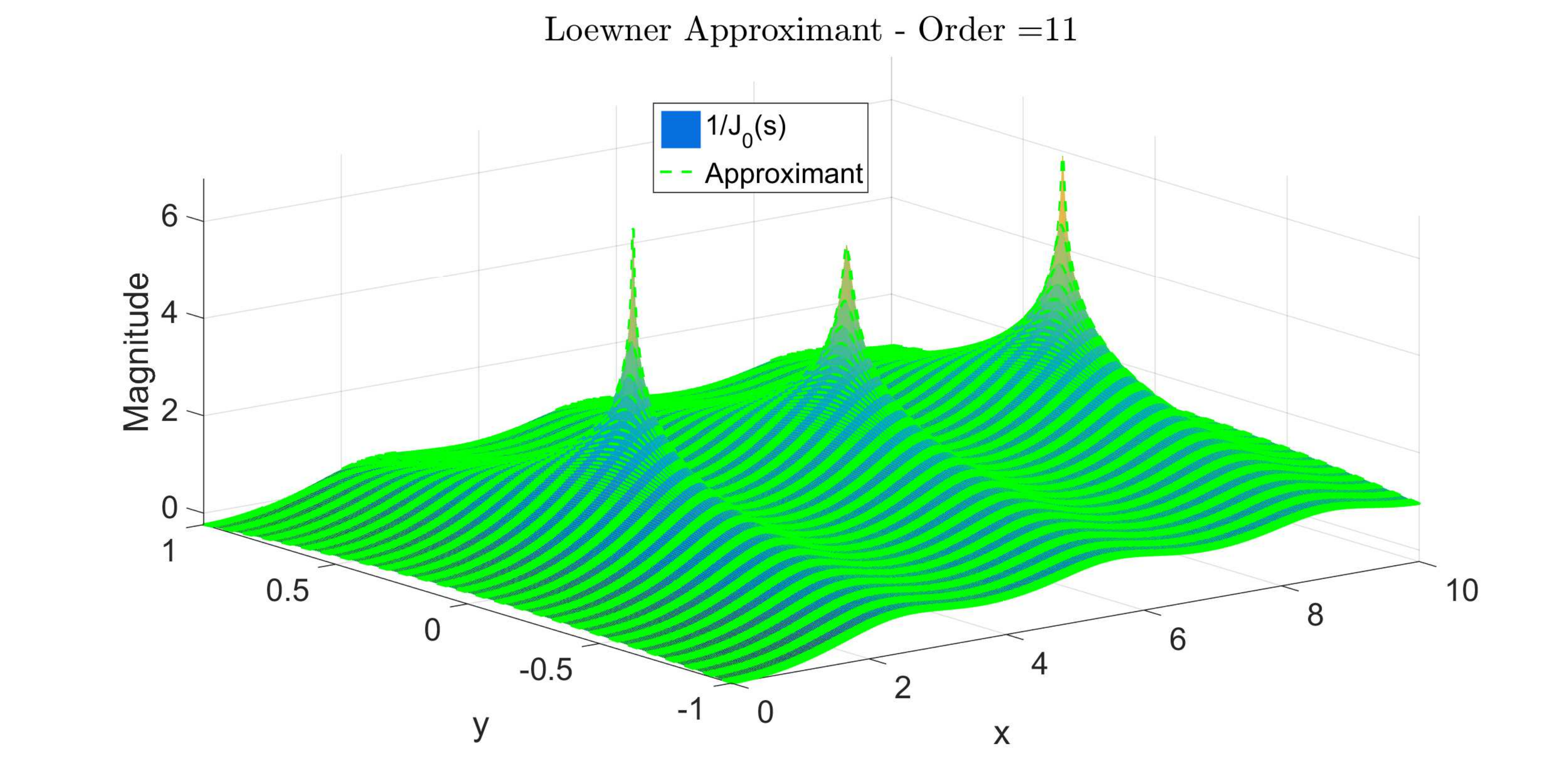}
  \caption{$H(s)$, with Loewner Approximant}
  \label{fig:sfig1}
\end{subfigure}%
\begin{subfigure}{.5\textwidth}
  \centering
  \includegraphics[width=1\linewidth]{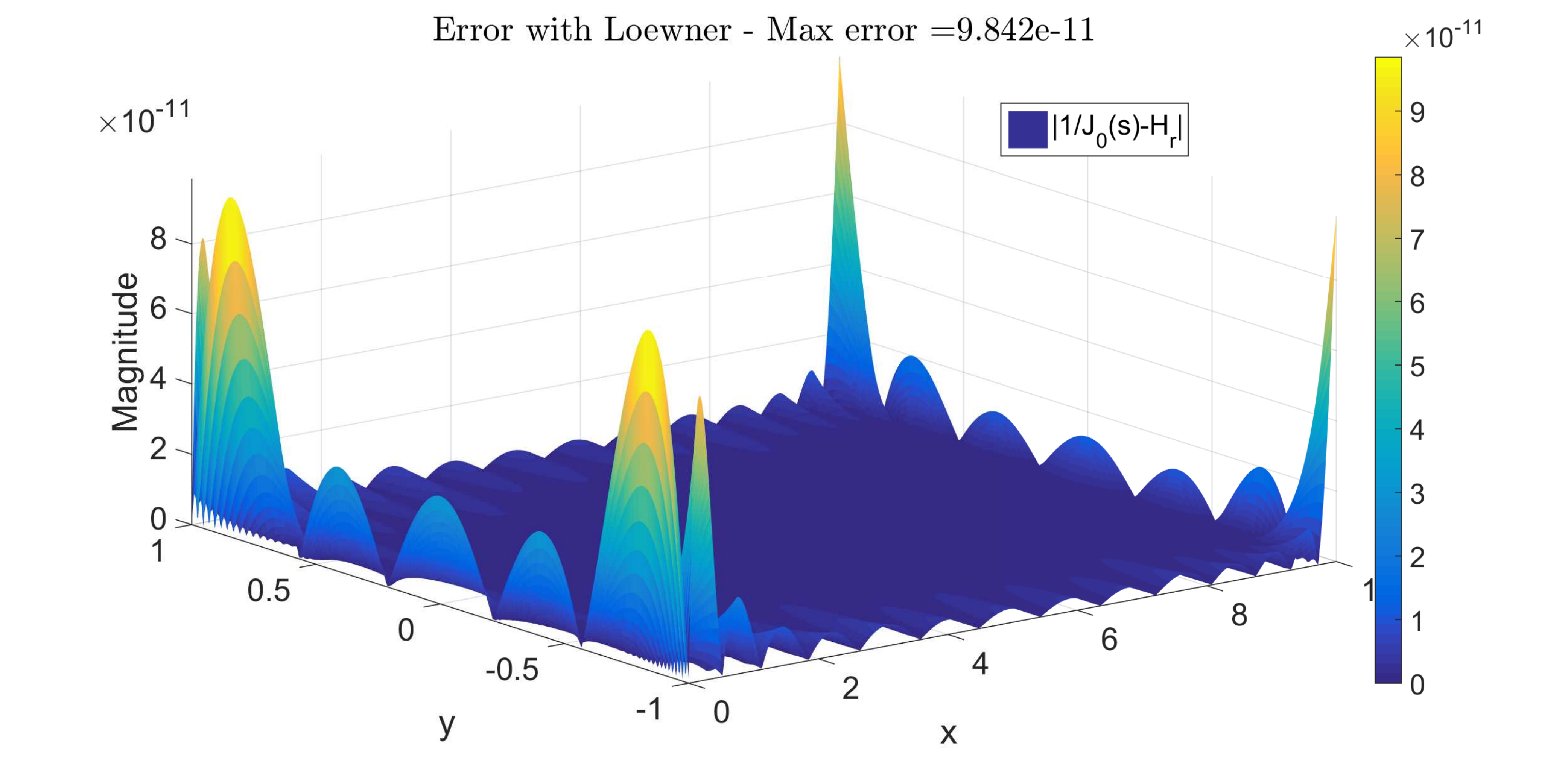}
  \caption{Absolute error of $|H(s)-H_{r}|$, over the $\bOmega$, domain.}
  \label{fig:sfig2}
\end{subfigure}
\caption{Approximant and the error over the dense domain $\bOmega_{grid}$.}
\end{figure}
\vspace{3mm}
\begin{figure}[h!]
\begin{subfigure}{.5\textwidth}
  \centering
  \includegraphics[width=1\linewidth]{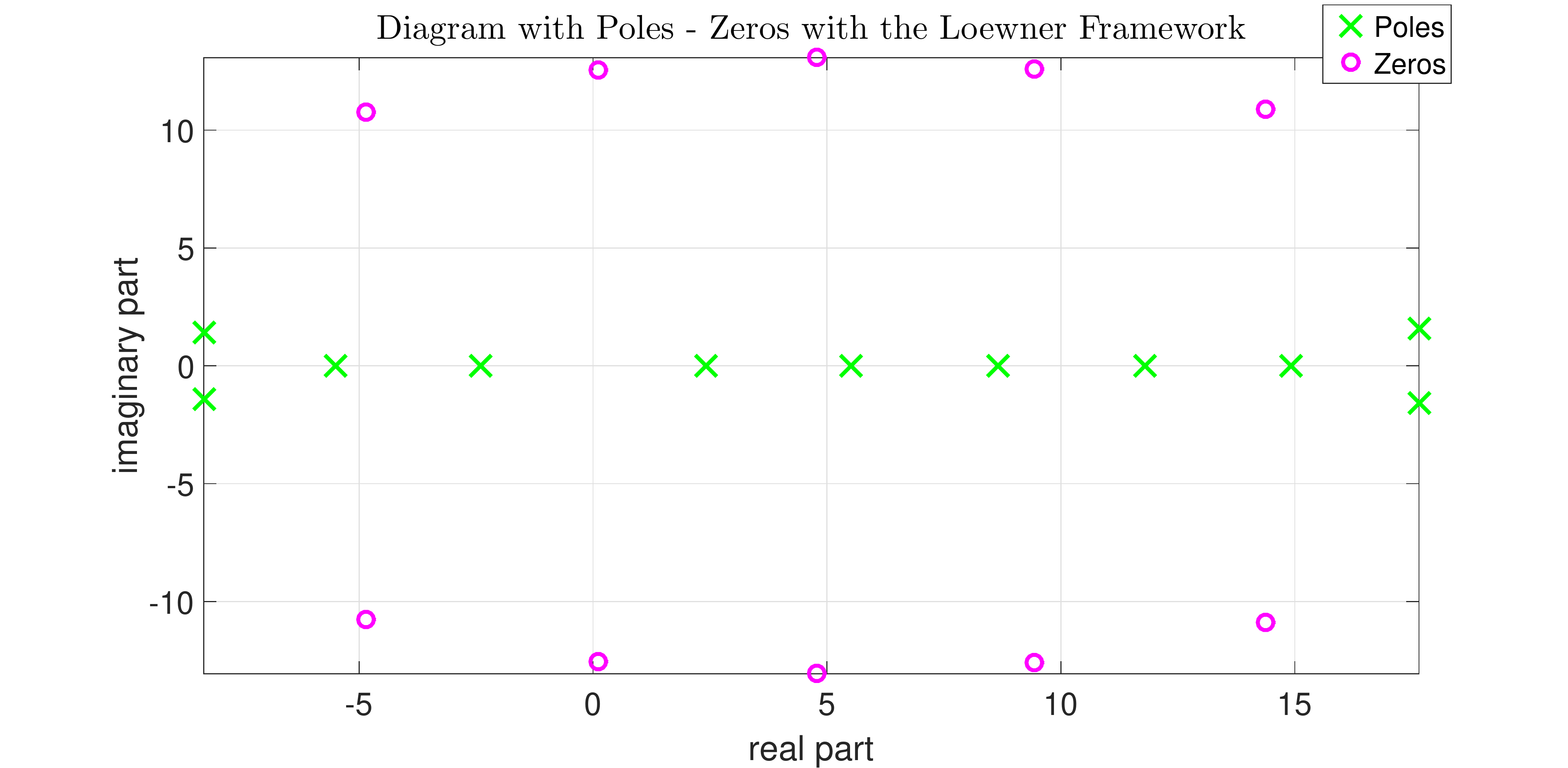}
  \caption{Poles $\color{green}X$, Zeros $\color{magenta}O$.}
  \label{fig:sfig1}
\end{subfigure}%
\begin{subfigure}{.5\textwidth}
  \centering
  \includegraphics[width=1\linewidth]{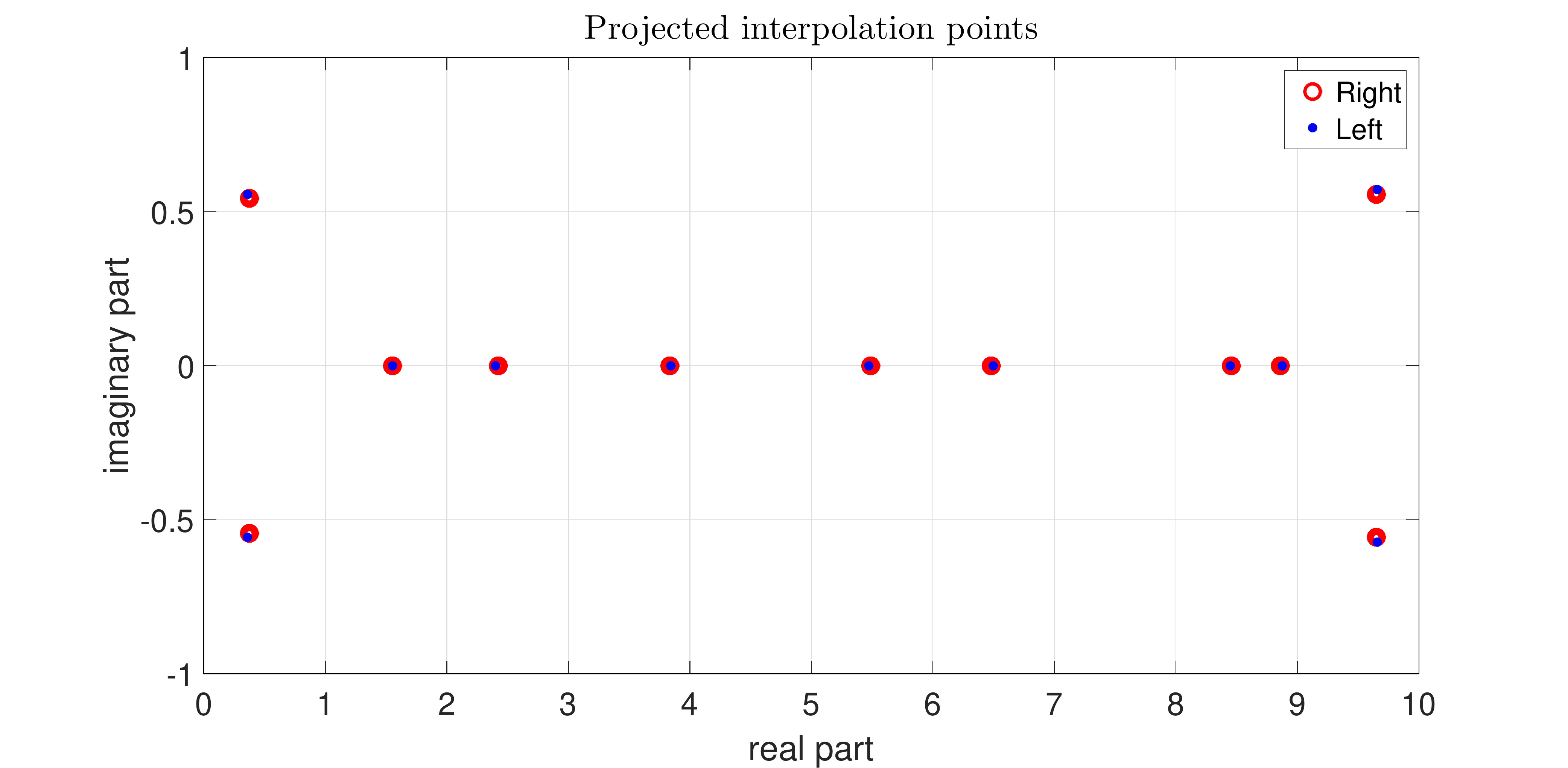}
  \caption{Right $\color{red}o$ and Left $\color{blue}\bullet$.}
  \label{fig:sfig2}
\end{subfigure}
\caption{Poles/zeros diagram and projected interpolation points in $\bOmega$.}
\end{figure}
Last, the corresponding poles/zeros and the projected left/right interpolation points are shown in the Figure 5. The real part of the poles is shown with 15 digits accuracy, while the imaginary part with 5. Also  the poles marked in bold  correspond to the original zeros of $J_{0}(s)$, in $\bOmega$, appendix 4. The ensuing left and right sets containing 11 projected interpolation points each (with 5 digits precision) are shown next.\\
\textbf{Remark: Thus, the initial 2121 interpolation points in $\bOmega$, are now compressed to 22}.

\begin{center}
\tiny
$\textbf{Poles}=\left(\begin{array}{c} -8.32213293322054- 1.4252\, \mathrm{i}\\ -8.32213289862456+ 1.4252\, \mathrm{i}\\ -5.51461491999547\\ -2.40481847965605\\ \textbf{2.40482555769577}\\ \textbf{5.52007811028631}\\ \textbf{8.65372791291101}\\ 11.7915356008908\\ 14.9135964357538\\ 17.6548692348549- 1.561\, \mathrm{i}\\ 17.654869354827+ 1.561\, \mathrm{i} \end{array}\right)$,
$\textbf{Zeros}=\left(\begin{array}{c} -4.8491 - 10.766\, \mathrm{i}\\ -4.8491 + 10.766\, \mathrm{i}\\ 0.12013 - 12.537\, \mathrm{i}\\ 0.12013 + 12.537\, \mathrm{i}\\ 4.785 + 13.066\, \mathrm{i}\\ 4.785 - 13.066\, \mathrm{i}\\ 9.4384 + 12.591\, \mathrm{i}\\ 9.4384 - 12.591\, \mathrm{i}\\ 14.367 + 10.868\, \mathrm{i}\\ 14.367 - 10.868\, \mathrm{i} \end{array}\right)$,
$\textbf{Right}=\left(\begin{array}{c} 0.37181 + 0.54303\, \mathrm{i}\\ 0.37181 - 0.54303\, \mathrm{i}\\ 1.5504 \,\\ 2.424 \,\\ 3.8331 \,\\ 5.4898 \\ 6.479 \\ 8.4519 \\ 8.8561 \\ 9.6519 - 0.55566\, \mathrm{i}\\ 9.6519 + 0.55566\, \mathrm{i} \end{array}\right)$,
$\textbf{Left}=\left(\begin{array}{c} 0.36164 - 0.55711\, \mathrm{i}\\ 0.36164 + 0.55711\, \mathrm{i}\\ 1.5491 \\ 2.4005 \\ 3.8418 \\ 5.473 \\ 6.4935 \\ 8.4516 \\ 8.8759 \\ 9.6605 - 0.57111\, \mathrm{i}\\ 9.6605 + 0.57111\, \mathrm{i} \end{array}\right)$.
\end{center}

\textbf{Trajectories of the projected interpolation points}. If we increase the number of interpolation points in $\bOmega$, the \textbf{left and right projected interpolation points seem to converge}. If we choose a specific way to partition the measurements (left and right), i.e. $\underset{i,j}{max}|\mu_{i}-\lambda_{j}|<\epsilon$, then we noticed that the corresponding pairs of projected points are $\epsilon$ order, close each other. The trajectories shown below use from 100 to almost 20000 initial interpolation points in $\bOmega$. By creating an increasingly dense $\bOmega_{grid}$, we observe the following convergence patterns. More specifically, we densify  symmetrically $\bOmega$ domain as follows: $a\times a, 2a\times 2a,...,na\times na$, with $a,n\in\mathbb{N}$, (i.e. $a=10$, and $n=14$) and for each grid we solve the set of equations (13). Only some pairs of trajectories out of 11 pairs are shown, as the remaining pairs have similar convergence behavior. Finally, it turns out that the projected interpolation points remain in $\bOmega$. 
\begin{figure}[h!]
\begin{subfigure}{.5\textwidth}
  \centering
  \includegraphics[width=1\linewidth]{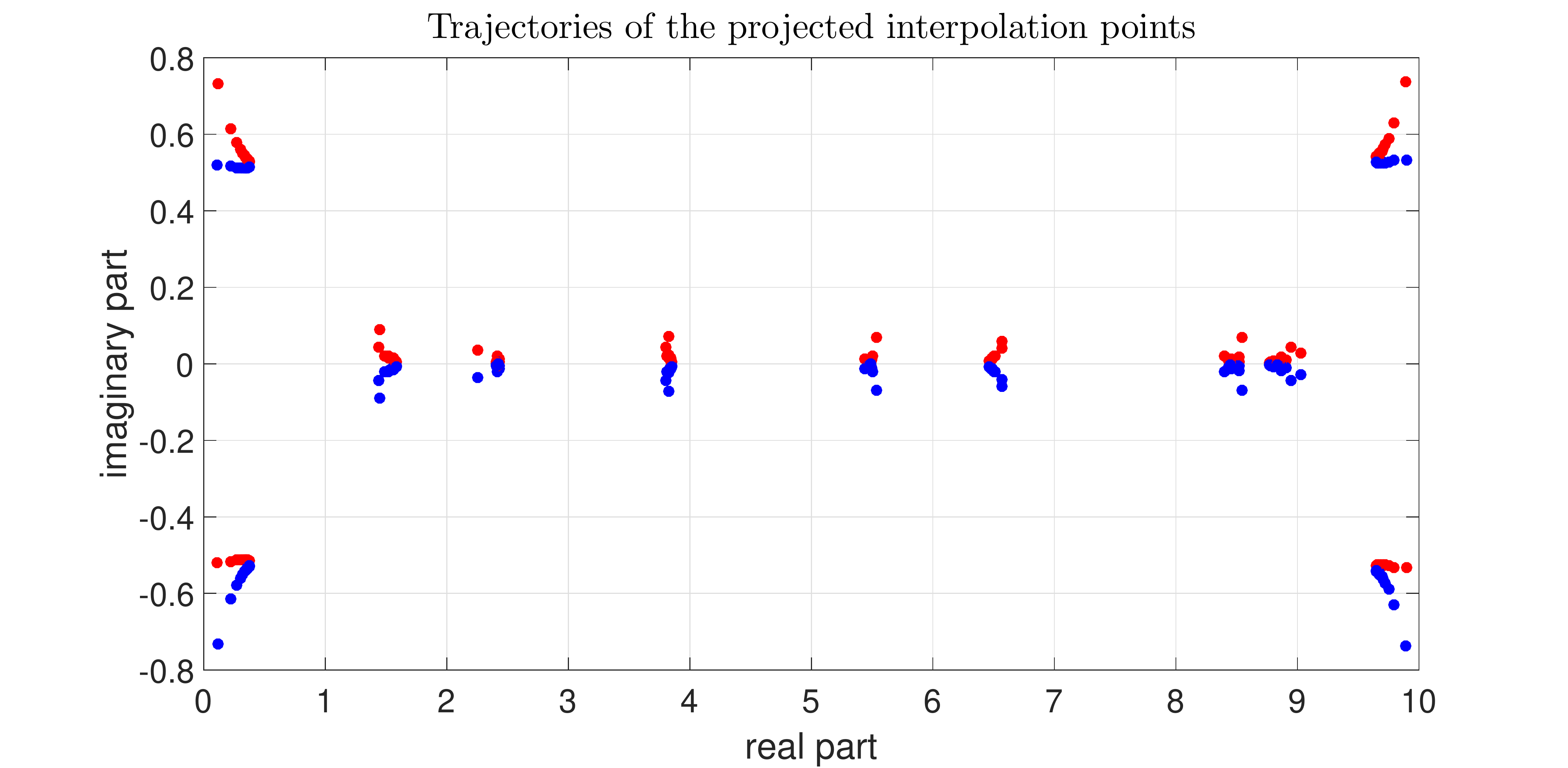}
  \caption{Trajectories of the projected interpolation points}
  \label{fig:sfig1}
\end{subfigure}%
\begin{subfigure}{.5\textwidth}
  \centering
  \includegraphics[width=1\linewidth]{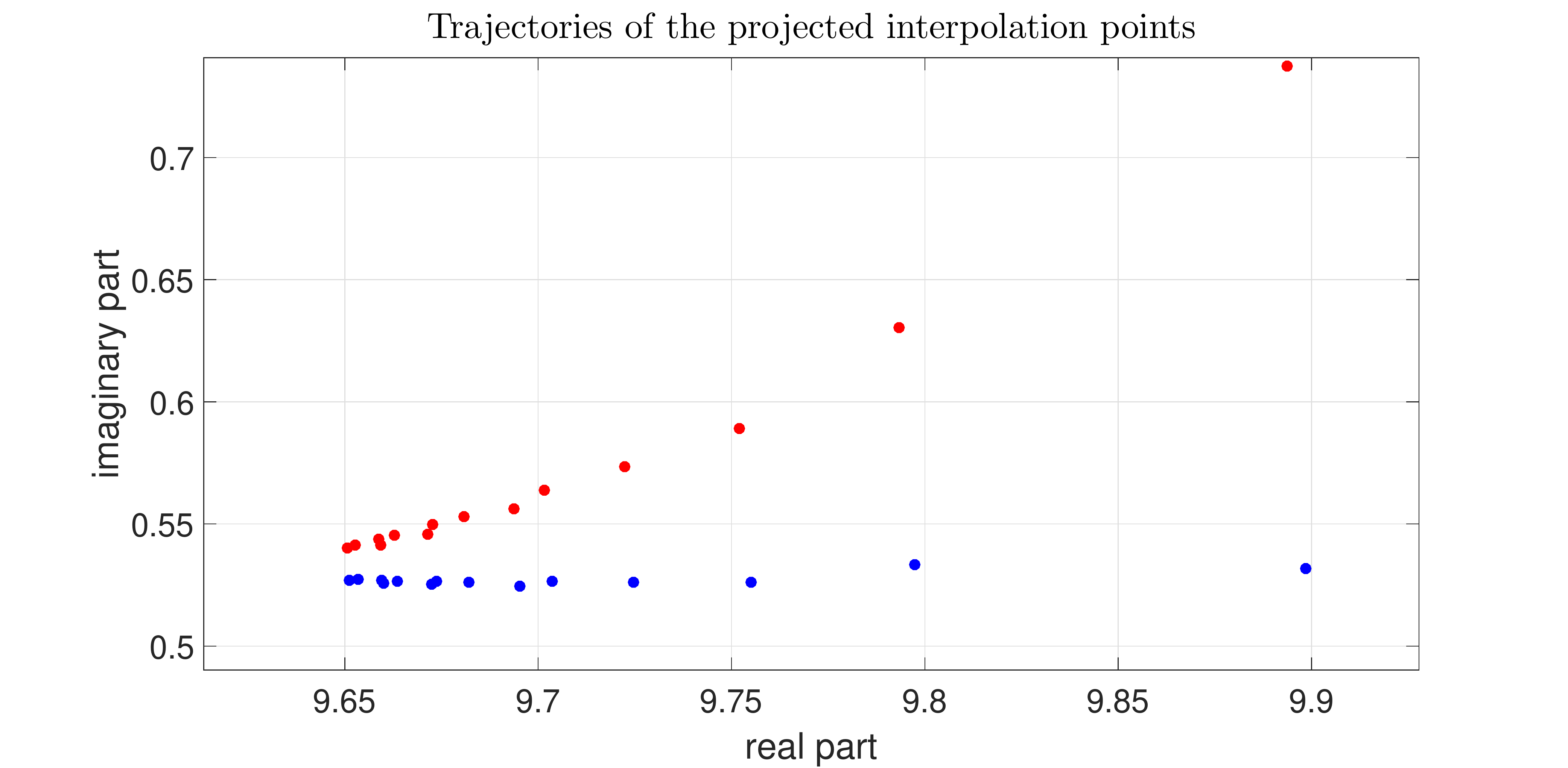}
  \caption{Zoom on the upper right of Figure 6a.}
  \label{fig:sfig2}
\end{subfigure}
\begin{subfigure}{.5\textwidth}
  \centering
  \includegraphics[width=1\linewidth]{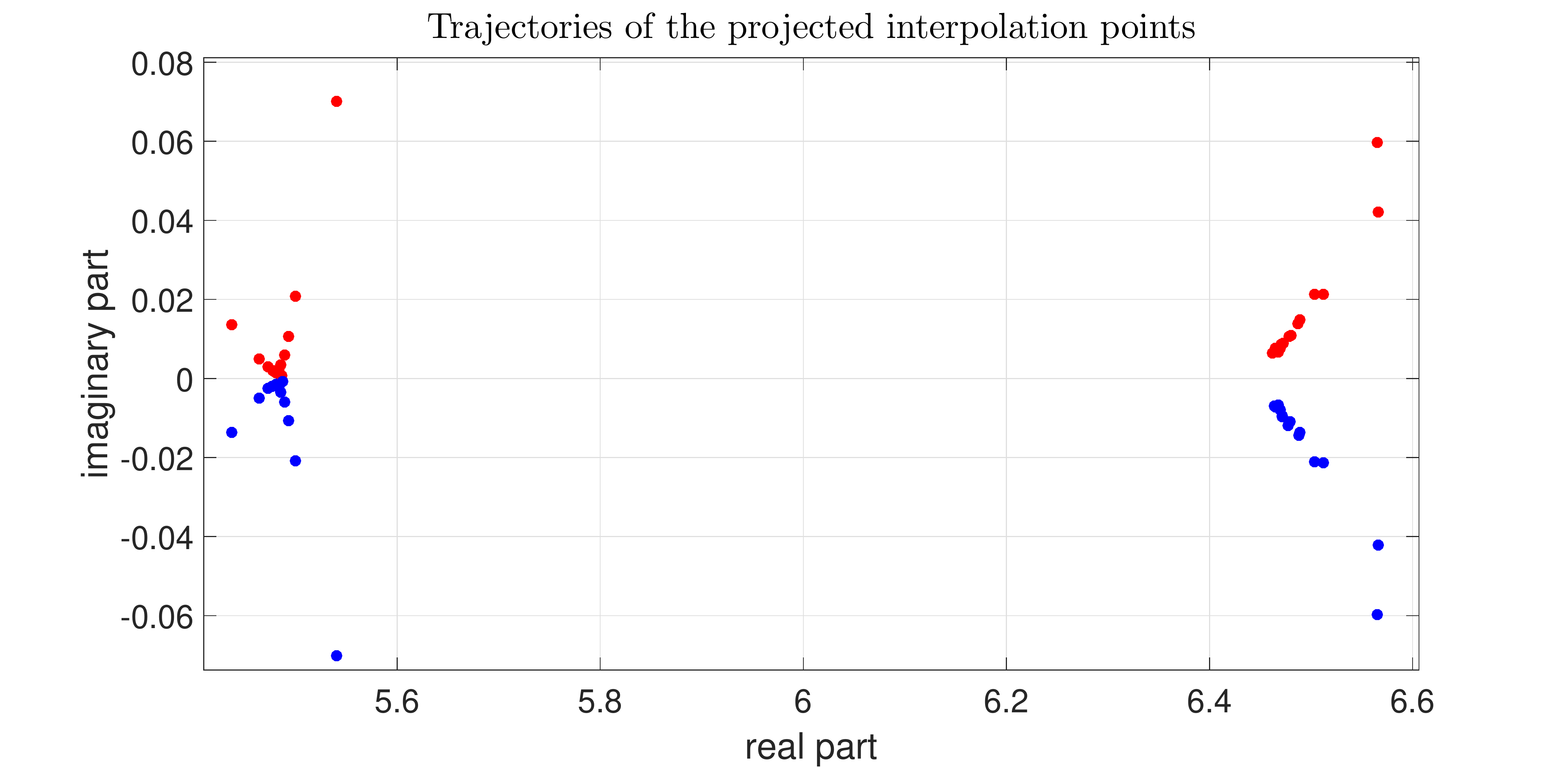}
  \caption{Convergence with real symmetry}
  \label{fig:sfig1}
\end{subfigure}%
\begin{subfigure}{.5\textwidth}
  \centering
  \includegraphics[width=1\linewidth]{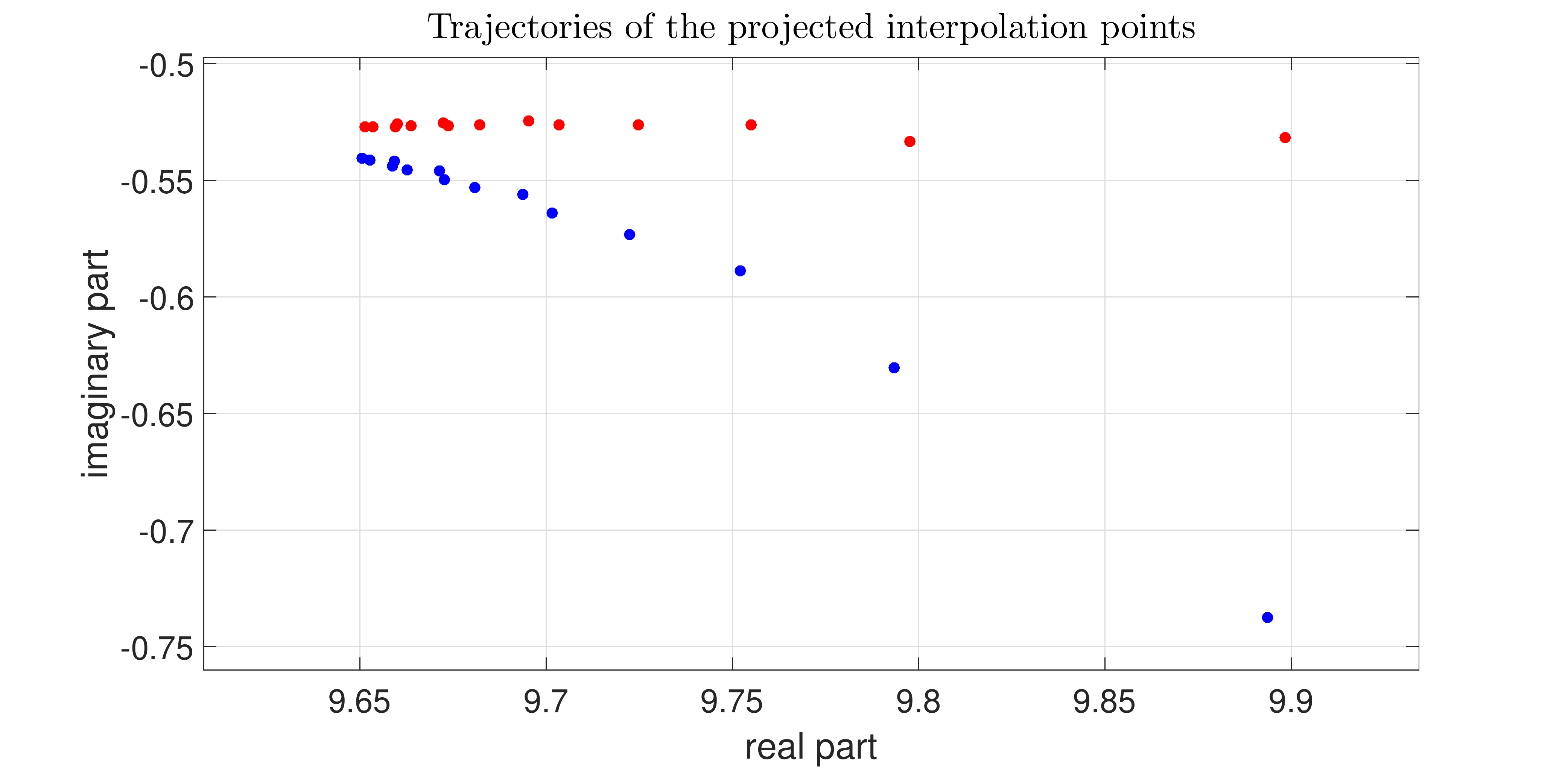}
  \caption{Conjugate trajectory of Figure 6b.}
  \label{fig:sfig1}
\end{subfigure}%
\caption{Pairs of trajectories of the projected interpolation points under increasingly denser grids. In some cases it seems that the convergence has oscillatory behavior. Also, left and right projected interpolation points seem to coincide when we select a particular way of separating the data (left-right).}
\label{fig:fig}
\end{figure}
\newpage
\subsubsection{The recursive Loewner algorithm to $H(s)$}
The measurements are as before. The resulting approximant compared with the original is as follows.
\begin{figure}[h!]
\begin{subfigure}{.5\textwidth}
  \centering
  \includegraphics[width=1\linewidth]{sampling2121pointsBessel-eps-converted-to.pdf}
  \caption{Corresponding sampling points of $H(s)$.}
  \label{fig:sfig1}
\end{subfigure}%
\begin{subfigure}{.5\textwidth}
  \centering
  \includegraphics[width=1\linewidth]{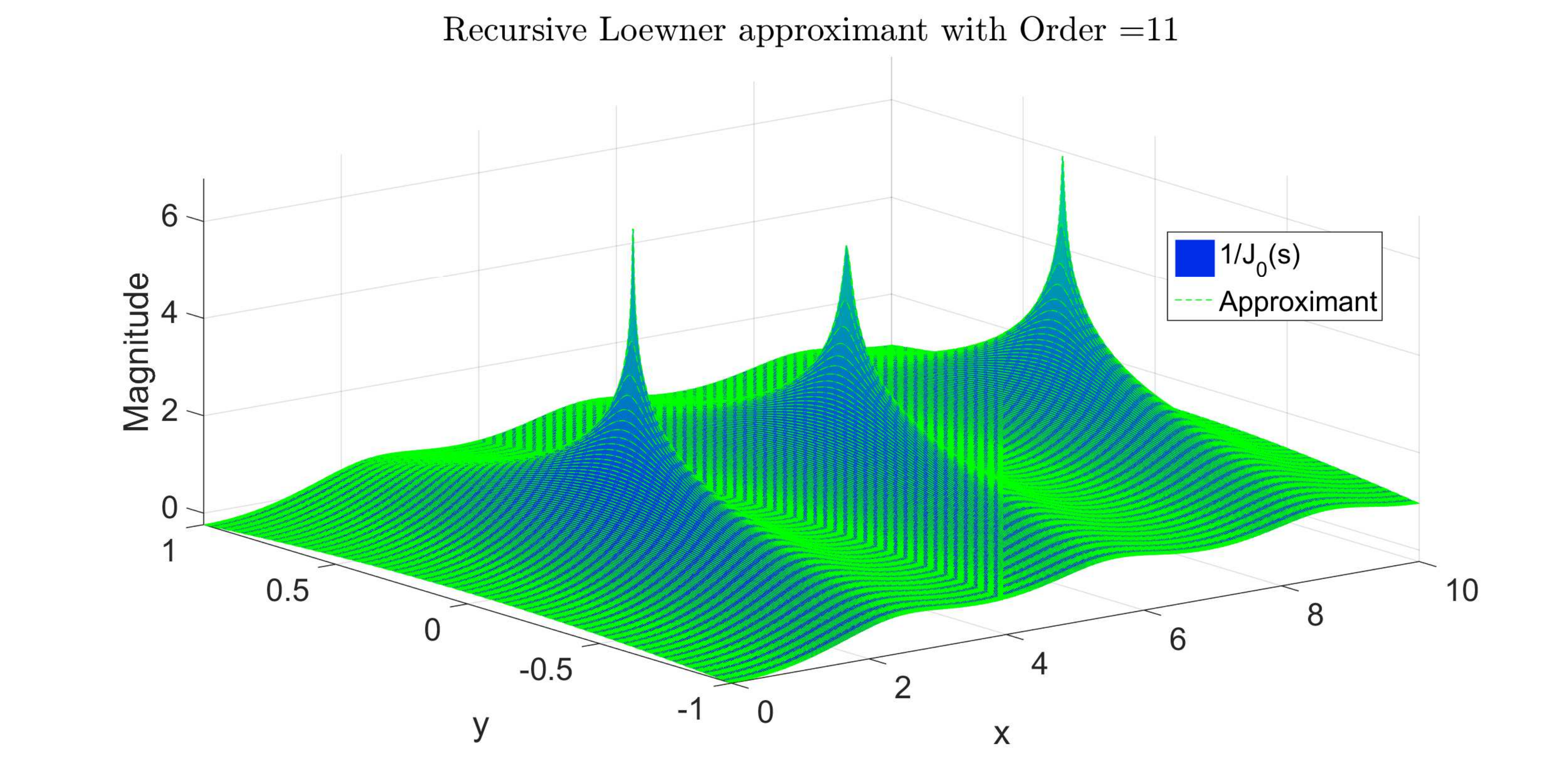}
  \caption{Approximant superimposed with $H(s)$, in $\bOmega$.}
  \label{fig:sfig2}
\end{subfigure}
\begin{subfigure}{.5\textwidth}
  \centering
  \includegraphics[width=1\linewidth]{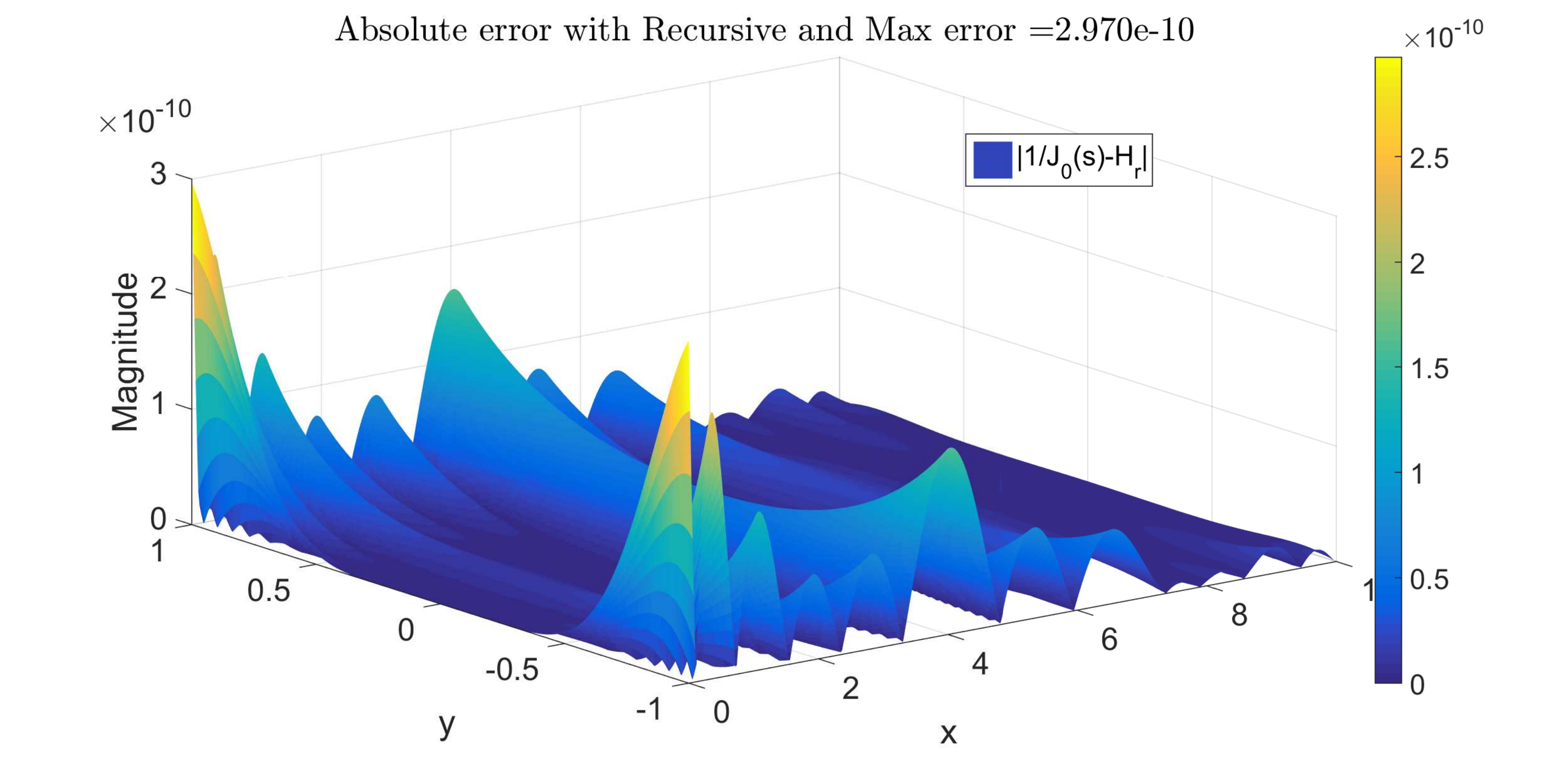}
  \caption{Error evaluation over the $\bOmega$.}
  \label{fig:sfig1}
\end{subfigure}%
\begin{subfigure}{.5\textwidth}
  \centering
  \includegraphics[width=1\linewidth]{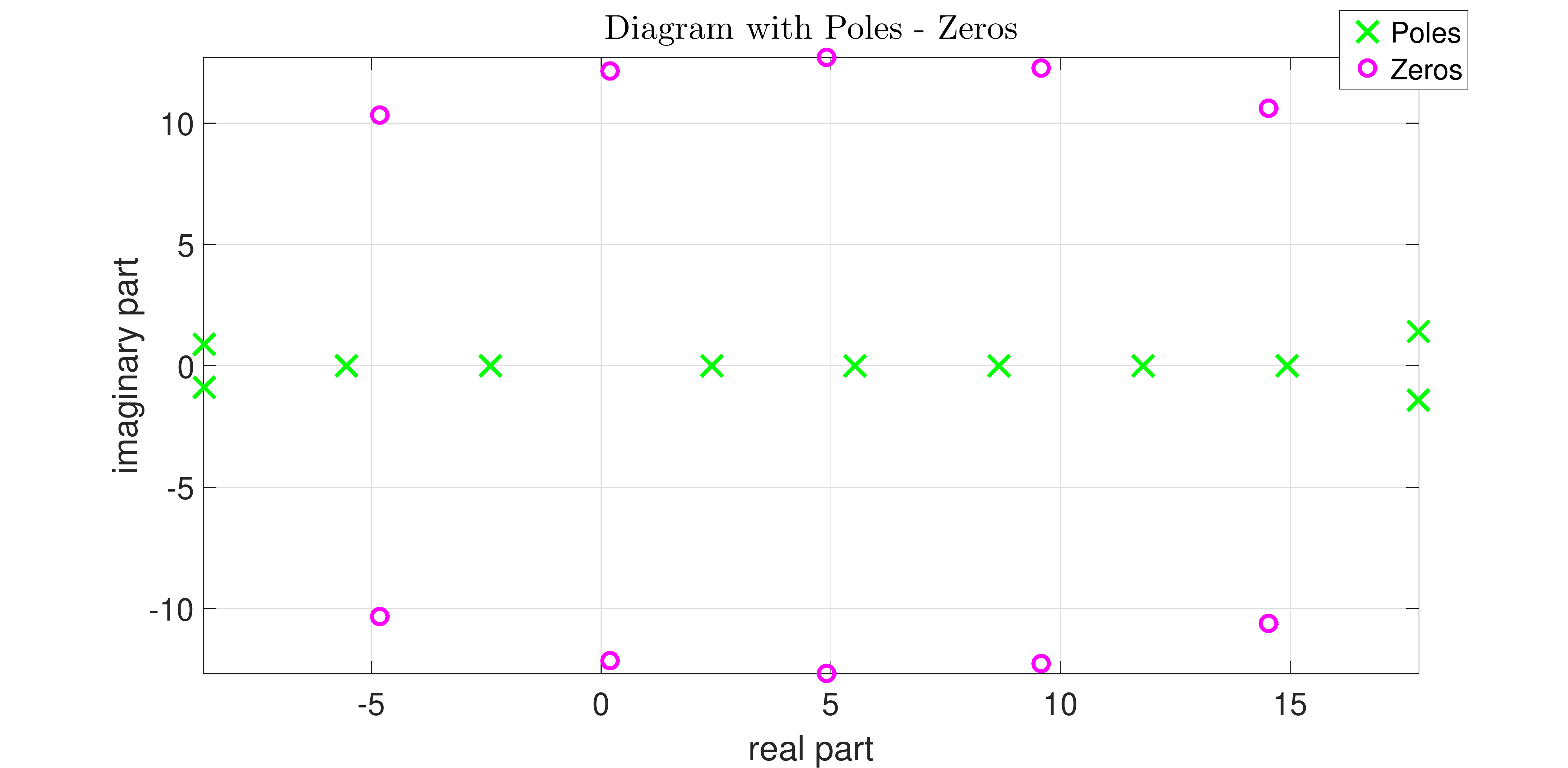}
  \caption{Poles and zeros diagram}
  \label{fig:sfig1}
\end{subfigure}%
\caption{Approximation results of the $H(s)$, over $\bOmega$, with the recursive Loewner approach.}
\label{fig:fig}
\end{figure}

In more detail the Figure 7a describes the interpolation points as shown in the left hand figure, section 2, page 5 with a depiction of these points on the corresponding samples of $H(s)$. The Figure 7b shows the 11th order approximant with the recursive Loewner superimposed on the plot of $H(s)$, over the dense grid $\bOmega_{grid}=[x_{1},...,x_{500}]\times[y_{1},...,y_{500}]\subset \bOmega$. As a result is the error plot in Figure 7c which shows for each point in the $\bOmega_{grid}$, the distance $|H_{r}(s)-H(s)|$, with $s \in \bOmega_{grid}$. The order of error is $O(10^{-10})$. Last, the corresponding poles/zeros are shown in the Figure 7d and for the exact values of them, someone can see Appendix 4 number 2.\\
\vspace{3mm}\\
\textbf{Remark 2.1.2.1:} With the recursive Loewner approach algorithm we need to choose the order of the approximant.\\
\vspace{1mm}\\
\textbf{Remark 2.1.2.2:} With the recursive Loewner approach algorithm usually we notice that the error is distributed more uniformly over $\bOmega$.  
\newpage
\subsubsection{The AAA algorithm applied to $H(s)$}
\begin{figure}[h]
\begin{subfigure}{.5\textwidth}
  \centering
  \includegraphics[width=1\linewidth]{sampling2121pointsBessel-eps-converted-to.pdf}
  \caption{Corresponding sampling points of $H(s)$.}
  \label{fig:sfig1}
\end{subfigure}%
\begin{subfigure}{.5\textwidth}
  \centering
  \includegraphics[width=1\linewidth]{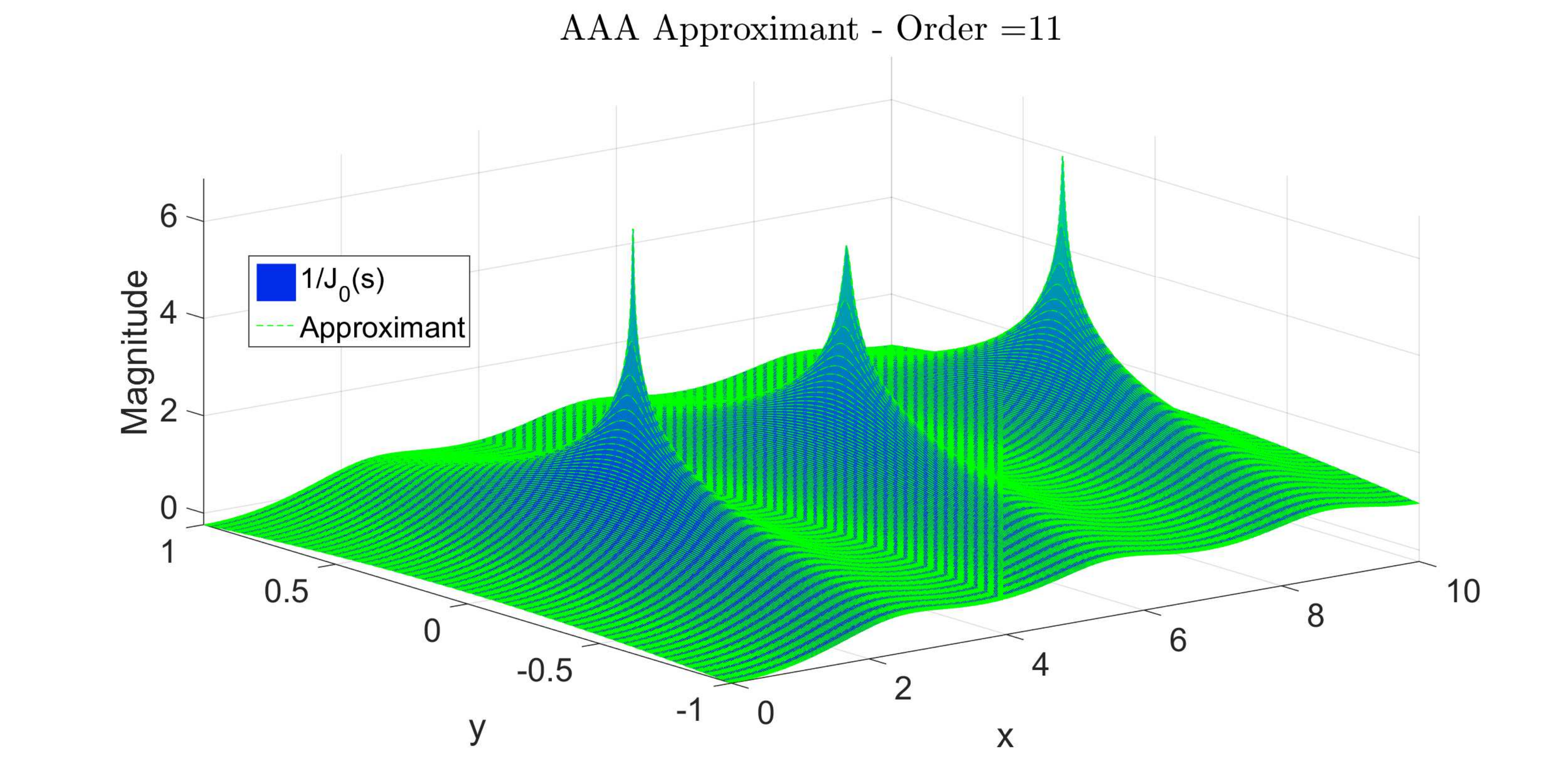}
  \caption{Approximant superimposed with $H(s)$, in $\bOmega$.}
  \label{fig:sfig2}
\end{subfigure}
\begin{subfigure}{.5\textwidth}
  \centering
  \includegraphics[width=1\linewidth]{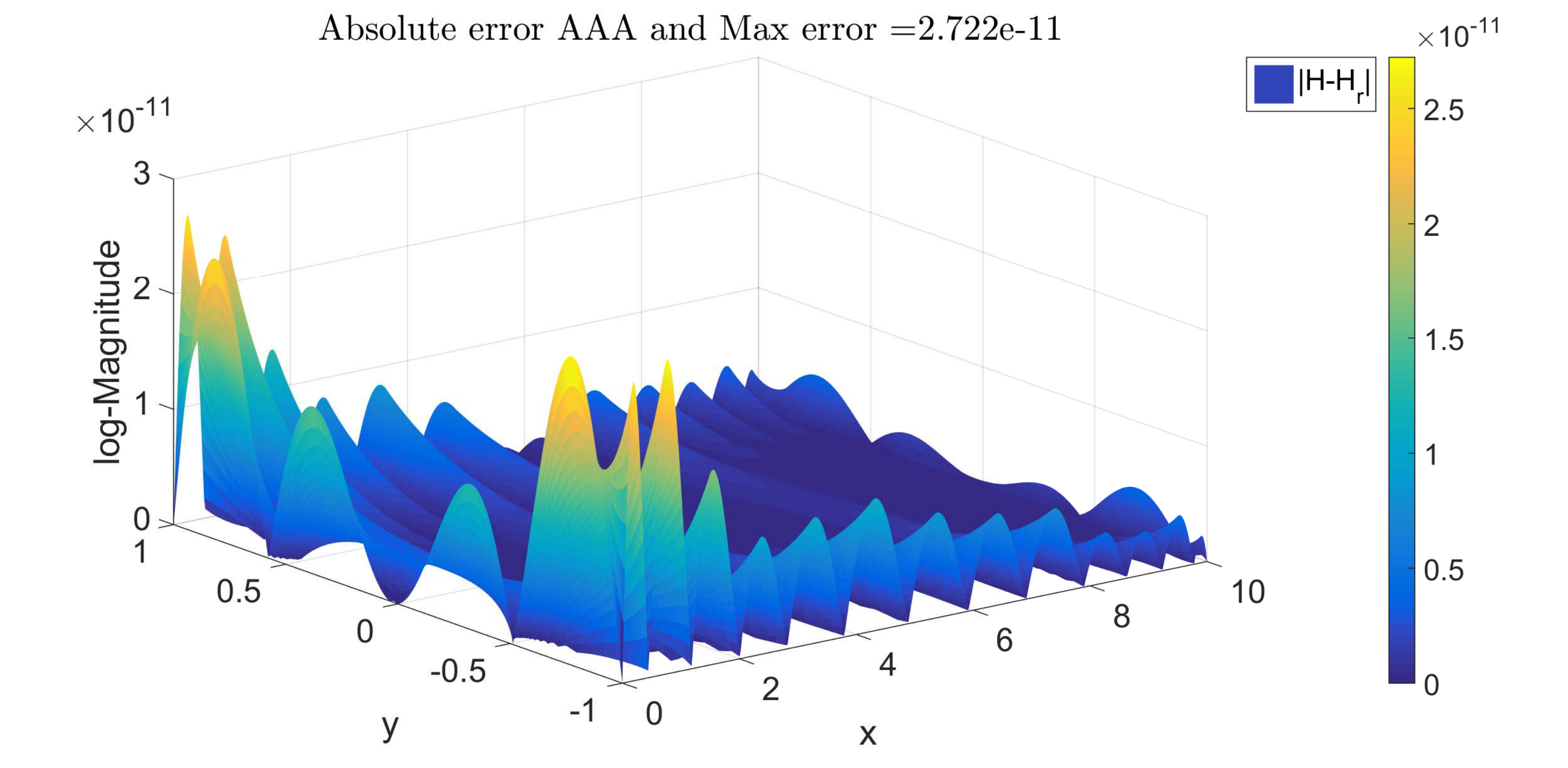}
  \caption{Error evaluation over the $\bOmega$ domain.}
  \label{fig:sfig1}
\end{subfigure}%
\begin{subfigure}{.5\textwidth}
  \centering
  \includegraphics[width=1\linewidth]{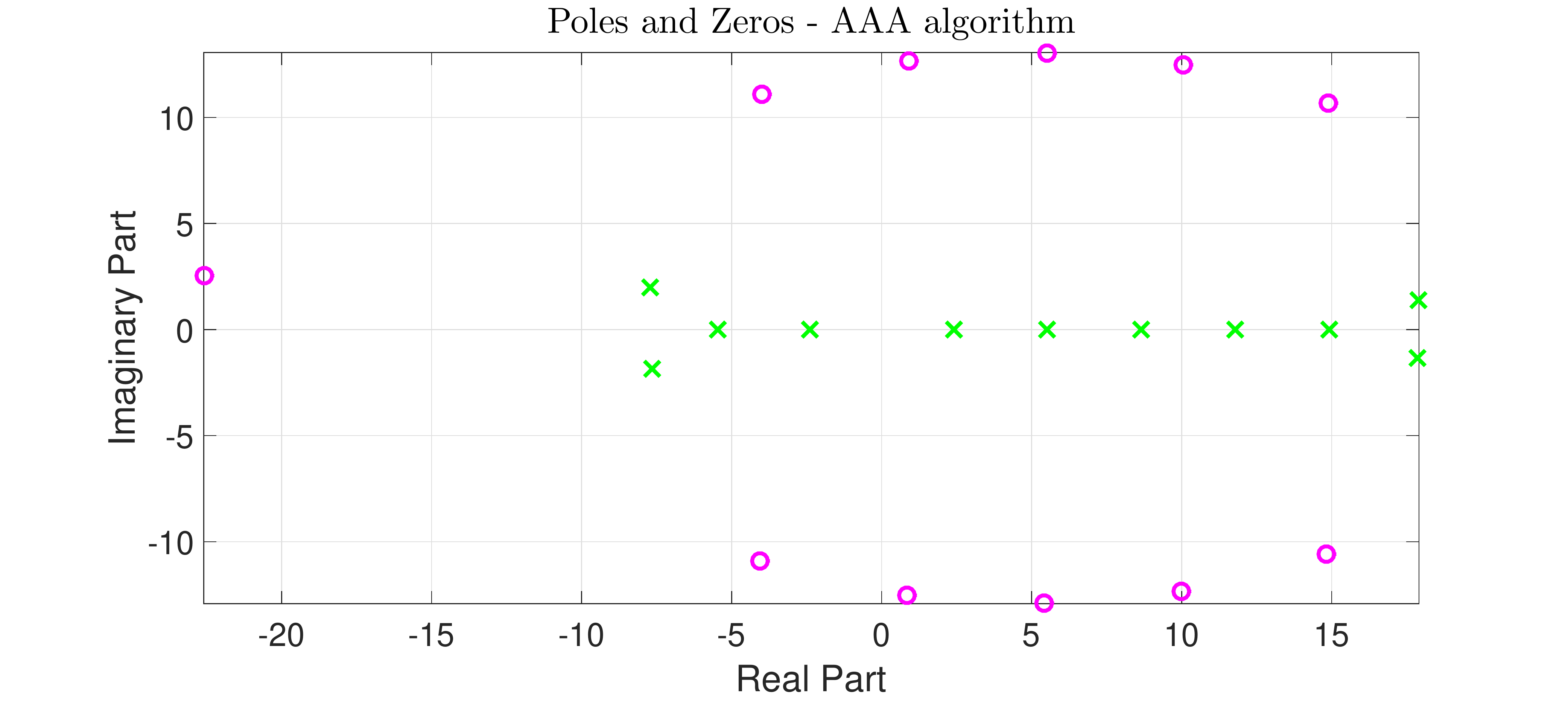}
  \caption{Poles and zeros diagram with AAA}
  \label{fig:sfig1}
\end{subfigure}%
\caption{Approximation results of the $H(s)$, over $\bOmega$, with the AAA algorithm approach.}
\label{fig:fig}
\end{figure}
In more detail the Figure 8a describes the interpolation points as shown in the left hand figure, section 2, page 5 with a depiction of these points on the corresponding samples of $H(s)$. The Figure 8b shows the 11th order approximant with the AAA approach superimposed on the plot of $H(s)$, over the dense grid $\bOmega_{grid}=[x_{1},...,x_{500}]\times[y_{1},...,y_{500}]\subset \bOmega$. As a result is the error plot in Figure 8c which shows for each point in the $\bOmega_{grid}$, the absolute distance $|H_{r}(s)-H(s)|$, with $s \in \bOmega_{grid}$. The order of error is $O(10^{-11})$. Last, the corresponding poles/zeros are shown in the Figure 8d and for the exact values, someone can see Appendix 4 number 3.\\

\begin{figure}
\centering
\includegraphics[width=0.7\linewidth, height=4cm]{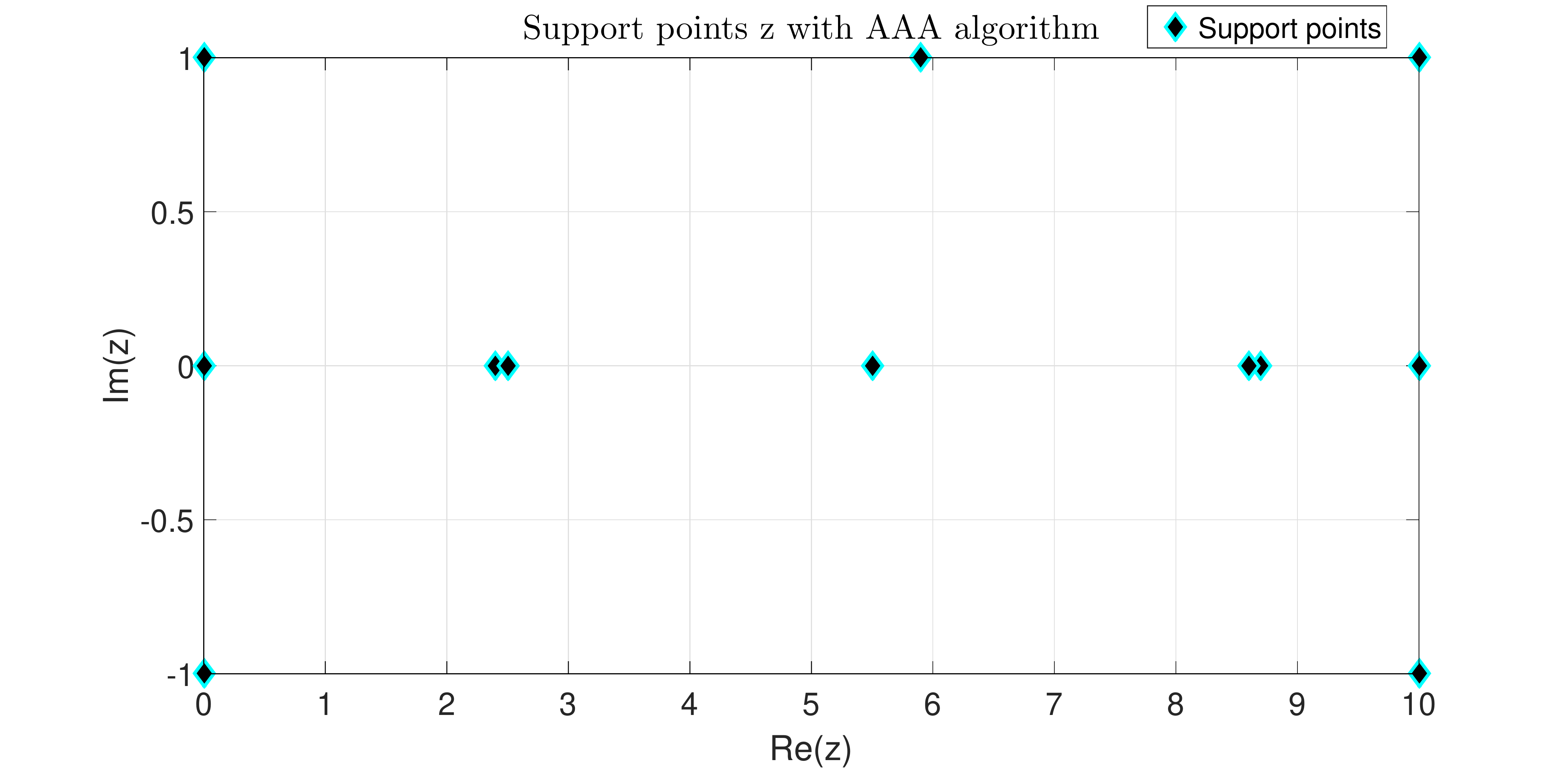} 
\caption{Support points from AAA algorithm in $\bOmega$.}
\end{figure}

\newpage
\subsubsection{The Vector Fitting (VF) method applied to $H(s)$}
The measurements are as before: 
\begin{figure}[h!]
\begin{subfigure}{.5\textwidth}
  \centering
  \includegraphics[width=1\linewidth]{sampling2121pointsBessel-eps-converted-to.pdf}
  \caption{Corresponding sampling points of $H(s)$.}
  \label{fig:sfig1}
\end{subfigure}%
\begin{subfigure}{.5\textwidth}
  \centering
  \includegraphics[width=1\linewidth]{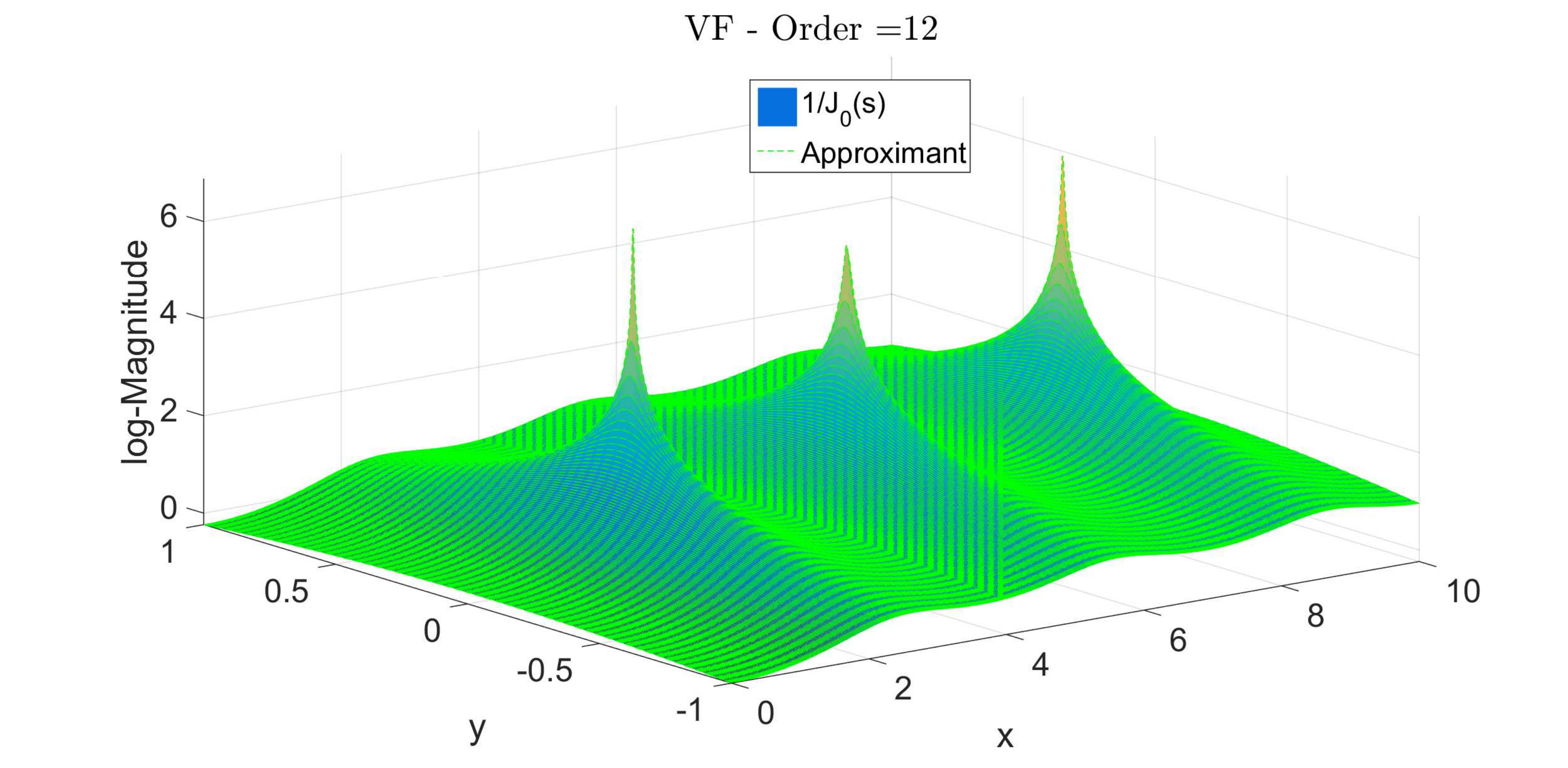}
  \caption{Approximant superimposed with $H(s)$, in $\bOmega$.}
  \label{fig:sfig2}
\end{subfigure}
\begin{subfigure}{.5\textwidth}
  \centering
  \includegraphics[width=1\linewidth]{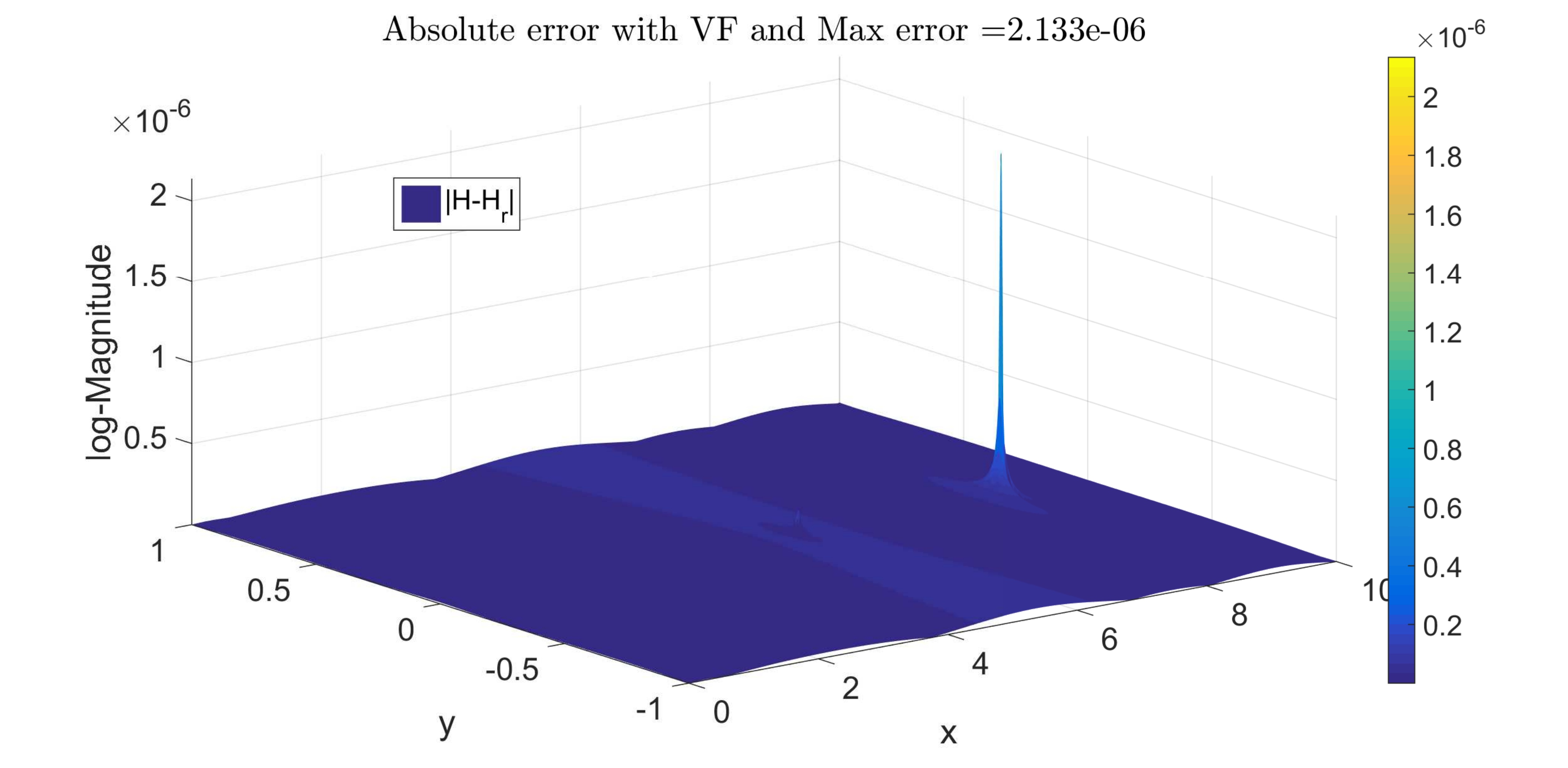}
  \caption{Error evaluation over the $\bOmega$.}
  \label{fig:sfig1}
\end{subfigure}%
\begin{subfigure}{.5\textwidth}
  \centering
  \includegraphics[width=1\linewidth]{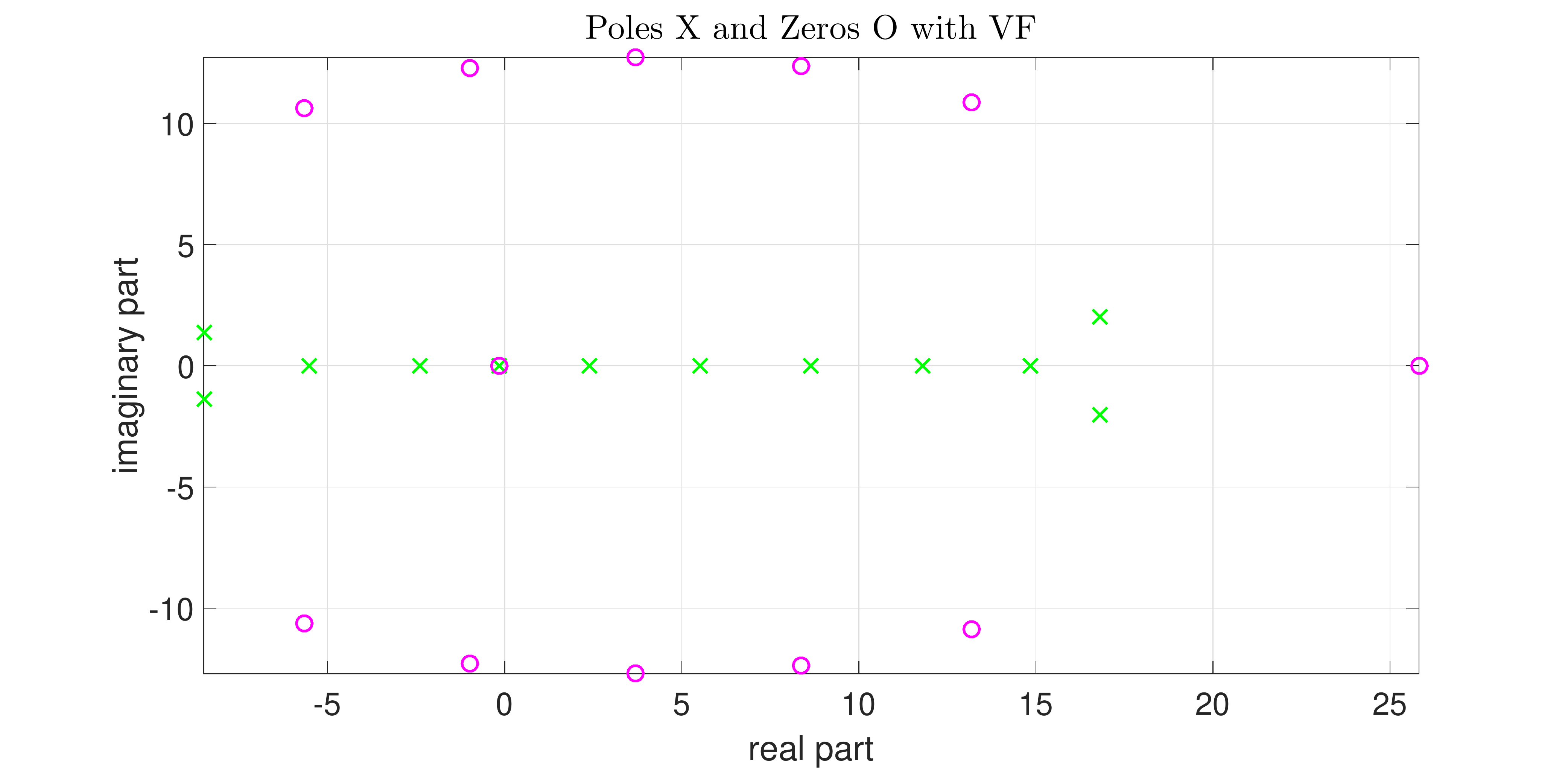}
  \caption{12 poles/zeros with one pole/zero cancellation}
  \label{fig:sfig1}
\end{subfigure}%
\caption{Approximation results of the $H(s)$, over $\bOmega$, with the Vector Fitting algorithm approach.}
\label{fig:fig}
\end{figure}

In more detail the Figure 10a describes the interpolation points as shown in the left hand figure, section 2, page 5 with a depiction of these points on the corresponding samples of $H(s)$. The Figure 10b shows the 11th order approximant with the VF approach superimposed on the plot of $H(s)$, over the dense grid $\bOmega_{grid}=[x_{1},...,x_{500}]\times[y_{1},...,y_{500}]\subset \bOmega$. As a result is the error plot in Figure 10c which shows for each point in the $\bOmega_{grid}$, the absolute distance $|H_{r}(s)-H(s)|$ with $s \in \bOmega_{grid}$. The order of error is $O(10^{-6})$. Last, the corresponding poles/zeros are shown in the Figure 10d and for the exact values, someone can see Appendix 4 number 4.\\
\vspace{2mm}\\
\textbf{Remark 2.1.4.1} Since the 5th pole and the first zero above are equal up to $10^{-9}$, (Figure 8d) or Appendix 4 number 4, they can be eliminated. After this \textbf{zero/pole} cancellation we obtain the same order of accuracy with an approximate of order $r=11$.\\
\vspace{2mm}\\
\textbf{Remark 2.1.4.2} We notice that the biggest error appears  in the vicinity of the 3rd pole. By changing the weights distribution we can improve the error locally. 
\\
\vspace{2mm}\\
\textbf{Remark 2.1.4.3} VF method also needs as an input the order of the approximant. 
\newpage
\subsection{Interpolation points as in figure 2b}
\subsubsection{The Loewner Framework applied to $H(s)$.}
We now make use of the uniformly distributed interpolation points as in Figure 2b.
\begin{figure}[h]
\begin{subfigure}{.5\textwidth}
  \centering
  \includegraphics[width=1\linewidth]{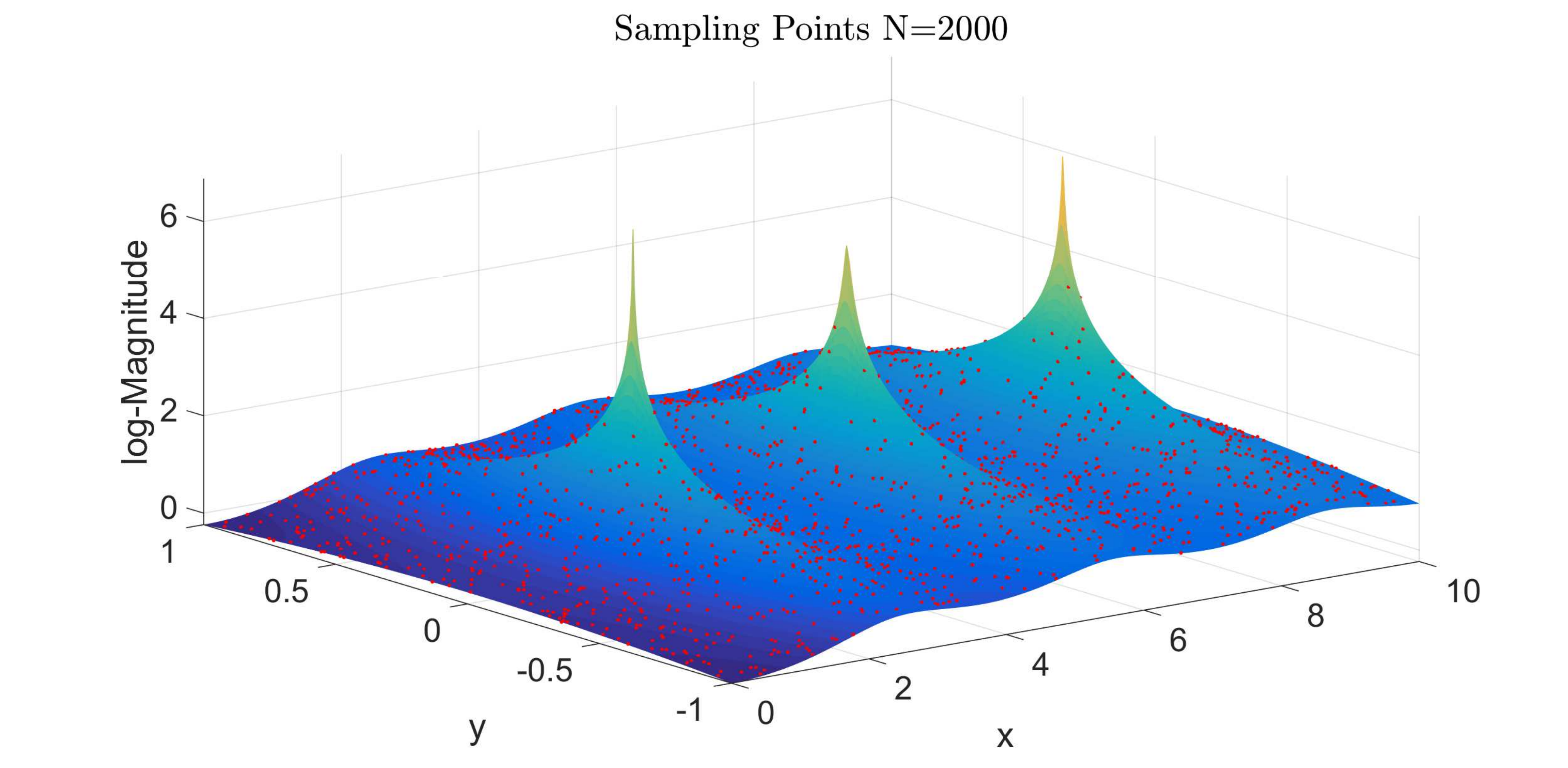}
  \caption{Corresponding sampling points of $H(s)$.}
  \label{fig:sfig1}
\end{subfigure}%
\begin{subfigure}{.5\textwidth}
  \centering
  \includegraphics[width=1\linewidth]{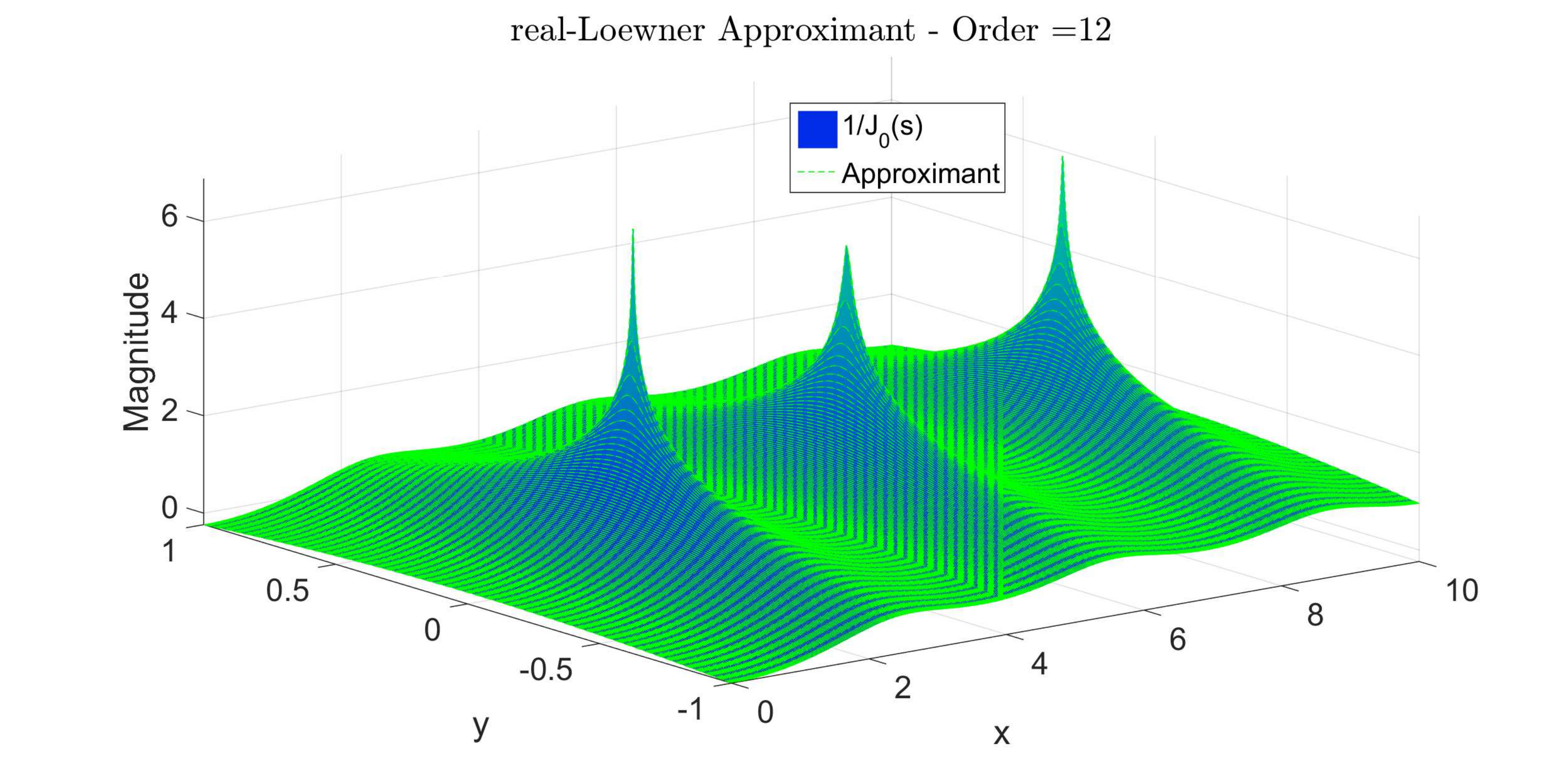}
  \caption{Approximant superimposed with $H(s)$, in $\bOmega$.}
  \label{fig:sfig2}
\end{subfigure}
\begin{subfigure}{.5\textwidth}
  \centering
  \includegraphics[width=1\linewidth]{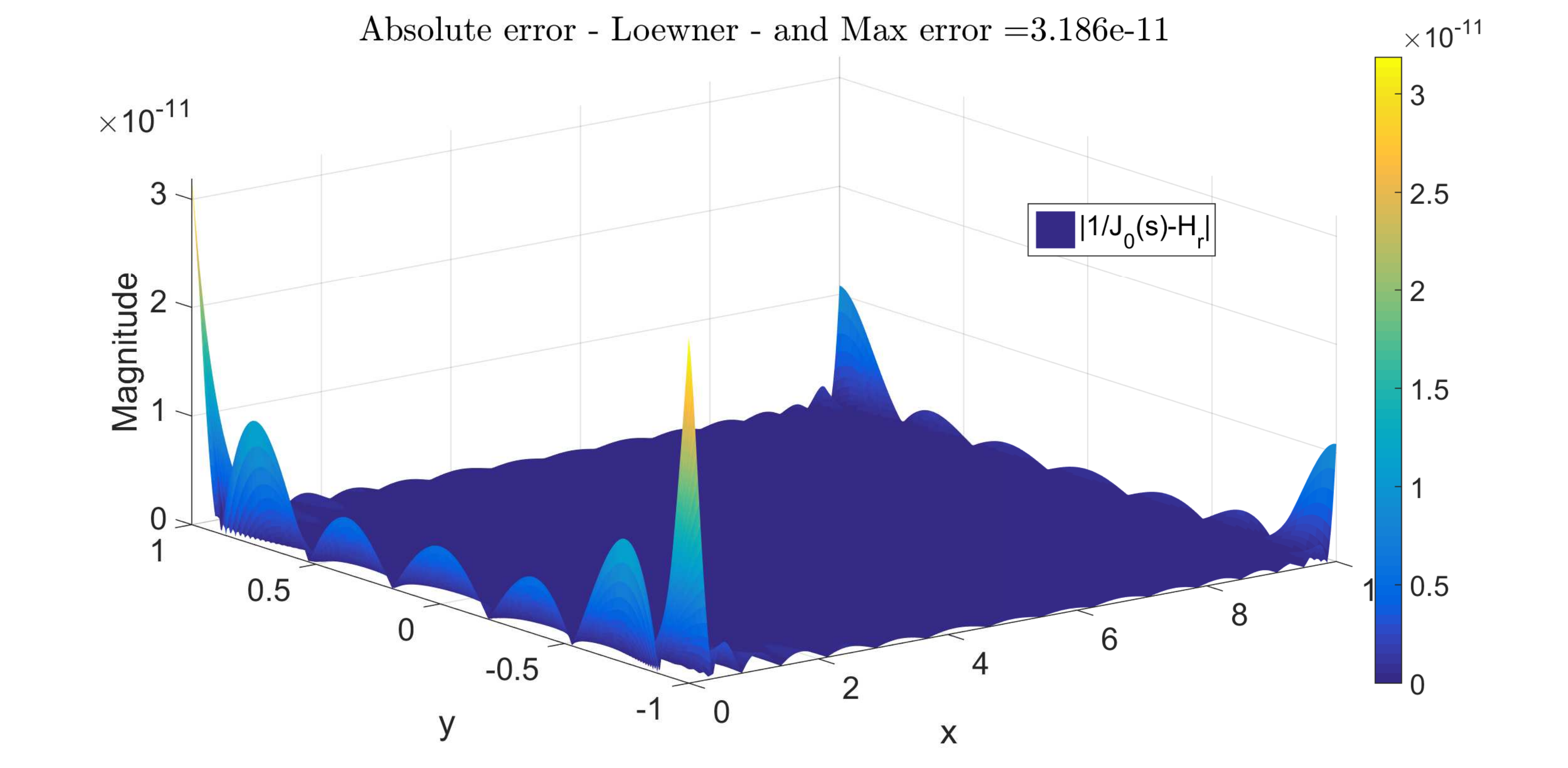}
  \caption{Error evaluation over the $\bOmega$.}
  \label{fig:sfig1}
\end{subfigure}%
\begin{subfigure}{.5\textwidth}
  \centering
  \includegraphics[width=1\linewidth]{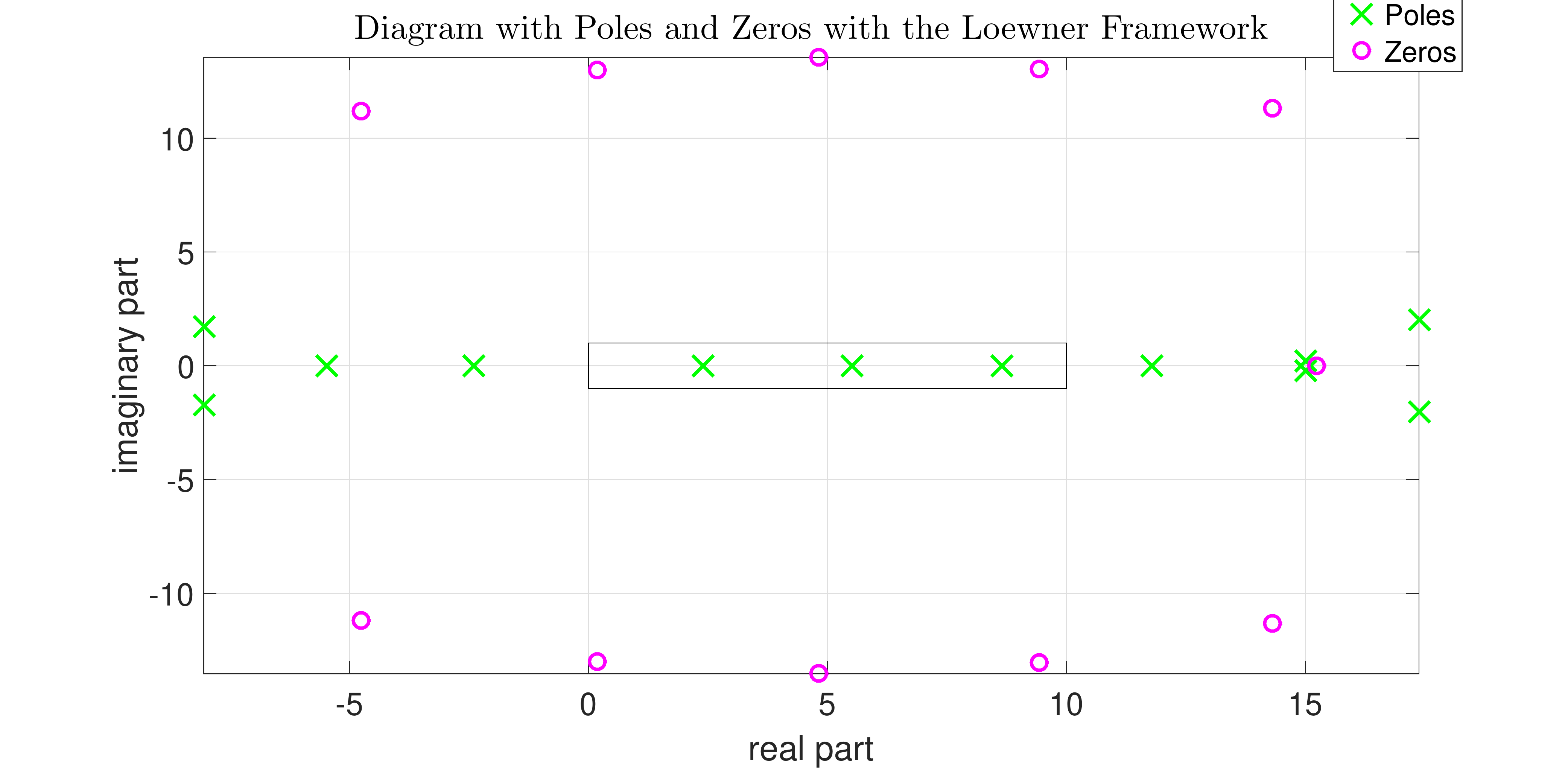}
  \caption{Poles and zeros diagram with $\bOmega$, rectangle}
  \label{fig:sfig1}
\end{subfigure}%
\caption{Approximation results of the $H(s)$, over $\bOmega$, with the Loewner framework approach.}
\label{fig:fig}
\end{figure}

In more detail the Figure 11a describes the interpolation points as shown in the right hand figure, section 2, page 5 with a depiction of these points on the corresponding samples of $H(s)$. The Figure 11b shows the 12th order approximant with the Loewner framework approach superimposed on the plot of $H(s)$, over the dense grid $\bOmega_{grid}=[x_{1},...,x_{500}]\times[y_{1},...,y_{500}]\subset \bOmega$. As a result we see the error plot in Figure 11c which shows for each point in the $\bOmega_{grid}$, the absolute distance $|H_{r}(s)-H(s)|$, with $s \in \bOmega_{grid}$. The order of error is $O(10^{-11})$. Last, the corresponding poles/zeros are shown in the Figure 11d and for the exact values, someone can see Appendix 4 number 5.\\
\begin{figure}[h]
\centering
\includegraphics[width=0.7\linewidth, height=4cm]{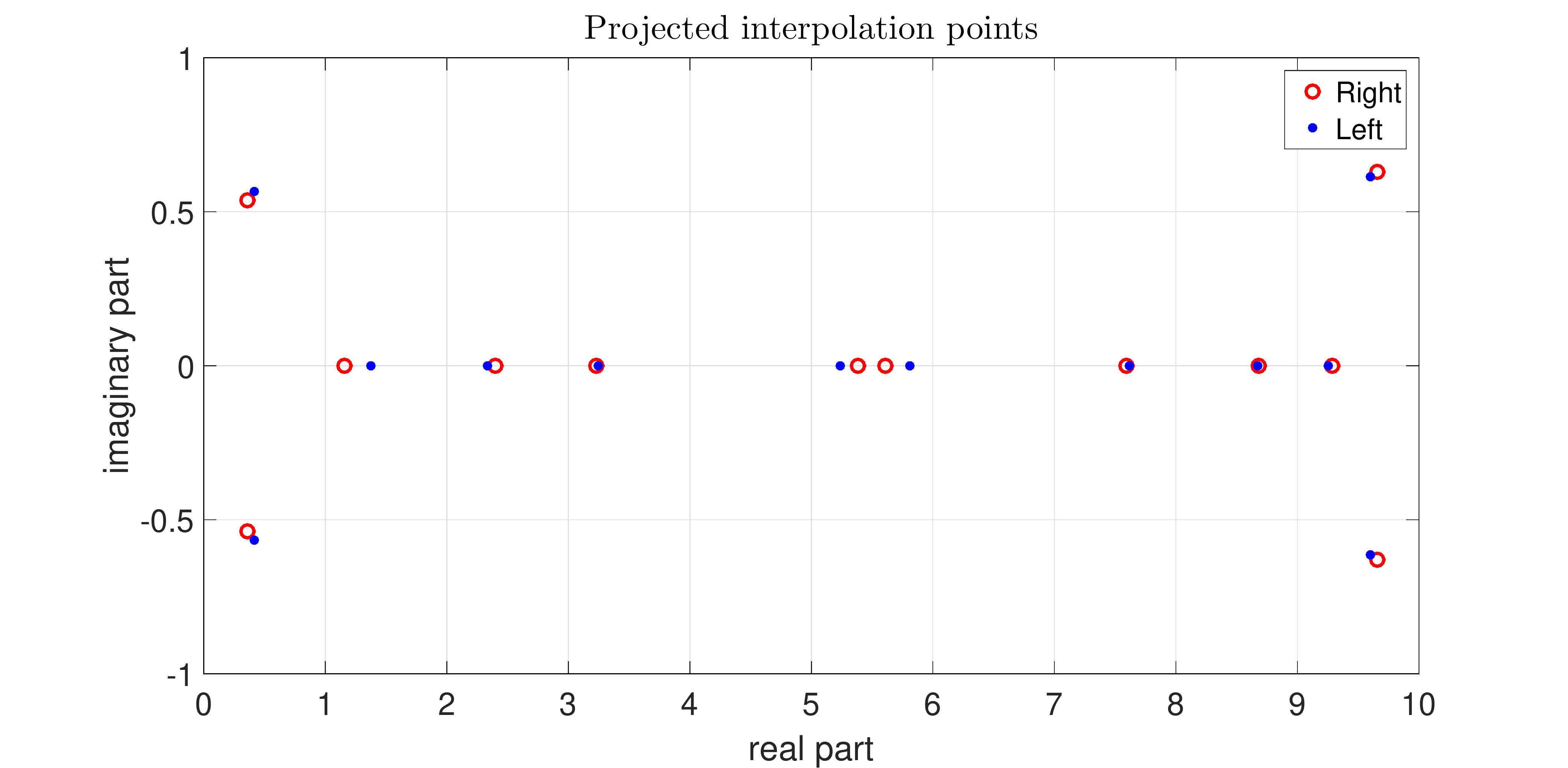} 
\caption{Left and right projected interpolation points in $\bOmega$.}
\end{figure}
\newpage
\subsubsection{The recursive Loewner Framework applied to $H(s)$.}
\begin{figure}[h!]
\begin{subfigure}{.5\textwidth}
  \centering
  \includegraphics[width=1\linewidth]{sampling2000pointsBessel-eps-converted-to.pdf}
  \caption{Corresponding sampling points of $H(s)$.}
  \label{fig:sfig1}
\end{subfigure}%
\begin{subfigure}{.5\textwidth}
  \centering
  \includegraphics[width=1\linewidth]{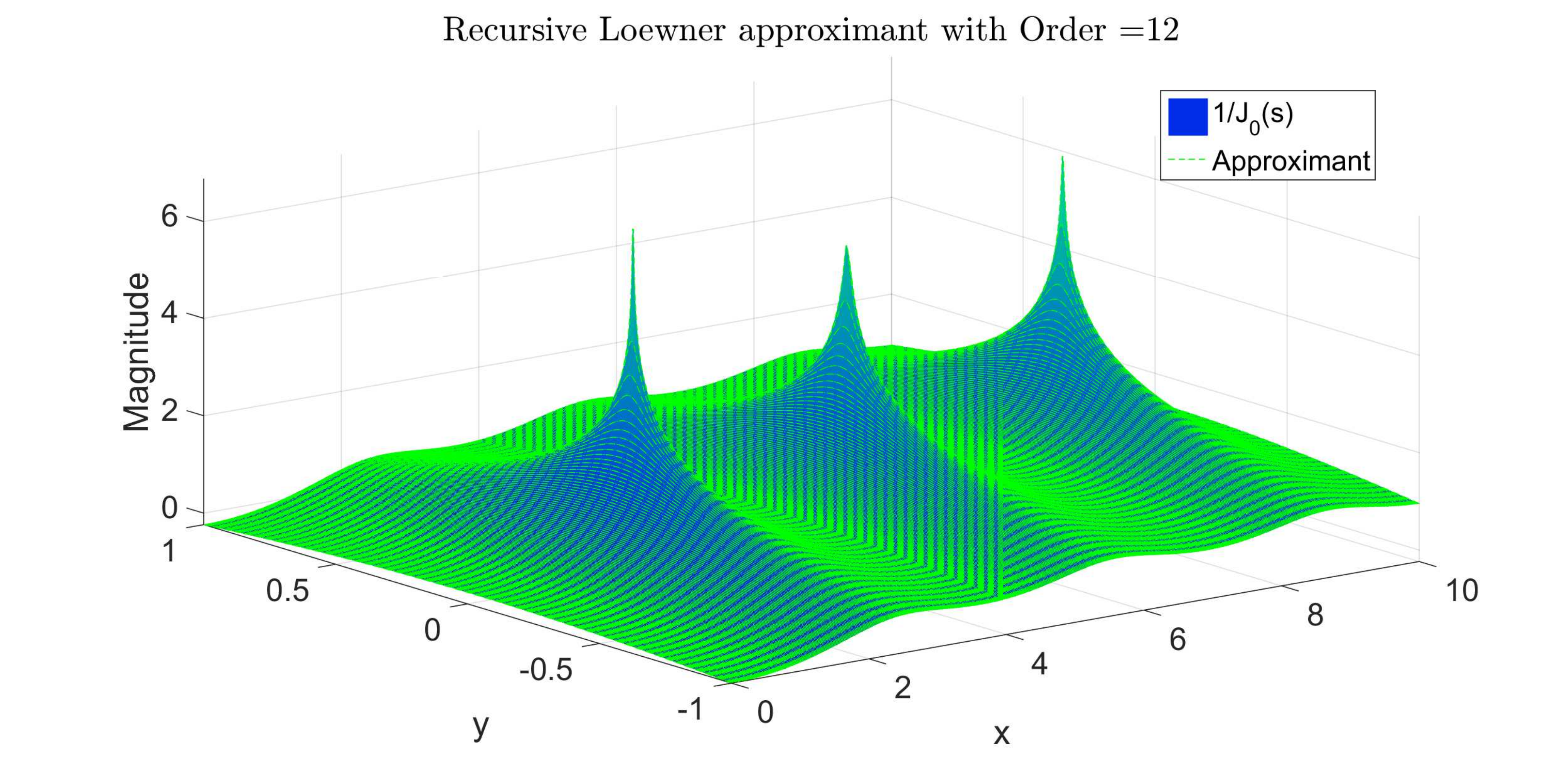}
  \caption{Approximant superimposed with $H(s)$, in $\bOmega$.}
  \label{fig:sfig2}
\end{subfigure}
\begin{subfigure}{.5\textwidth}
  \centering
  \includegraphics[width=1\linewidth]{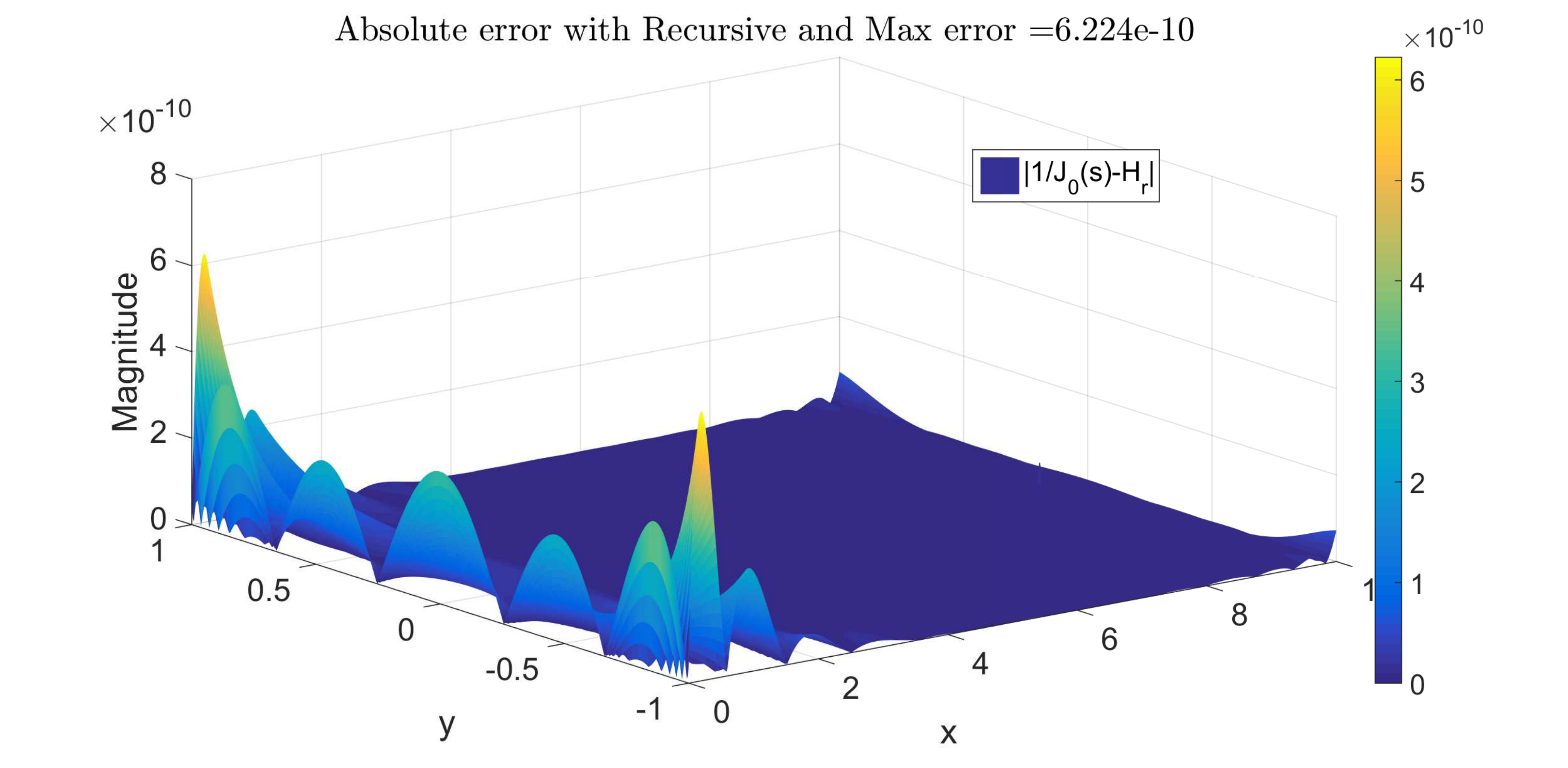}
  \caption{Error evaluation over the $\bOmega$.}
  \label{fig:sfig1}
\end{subfigure}%
\begin{subfigure}{.5\textwidth}
  \centering
  \includegraphics[width=1\linewidth]{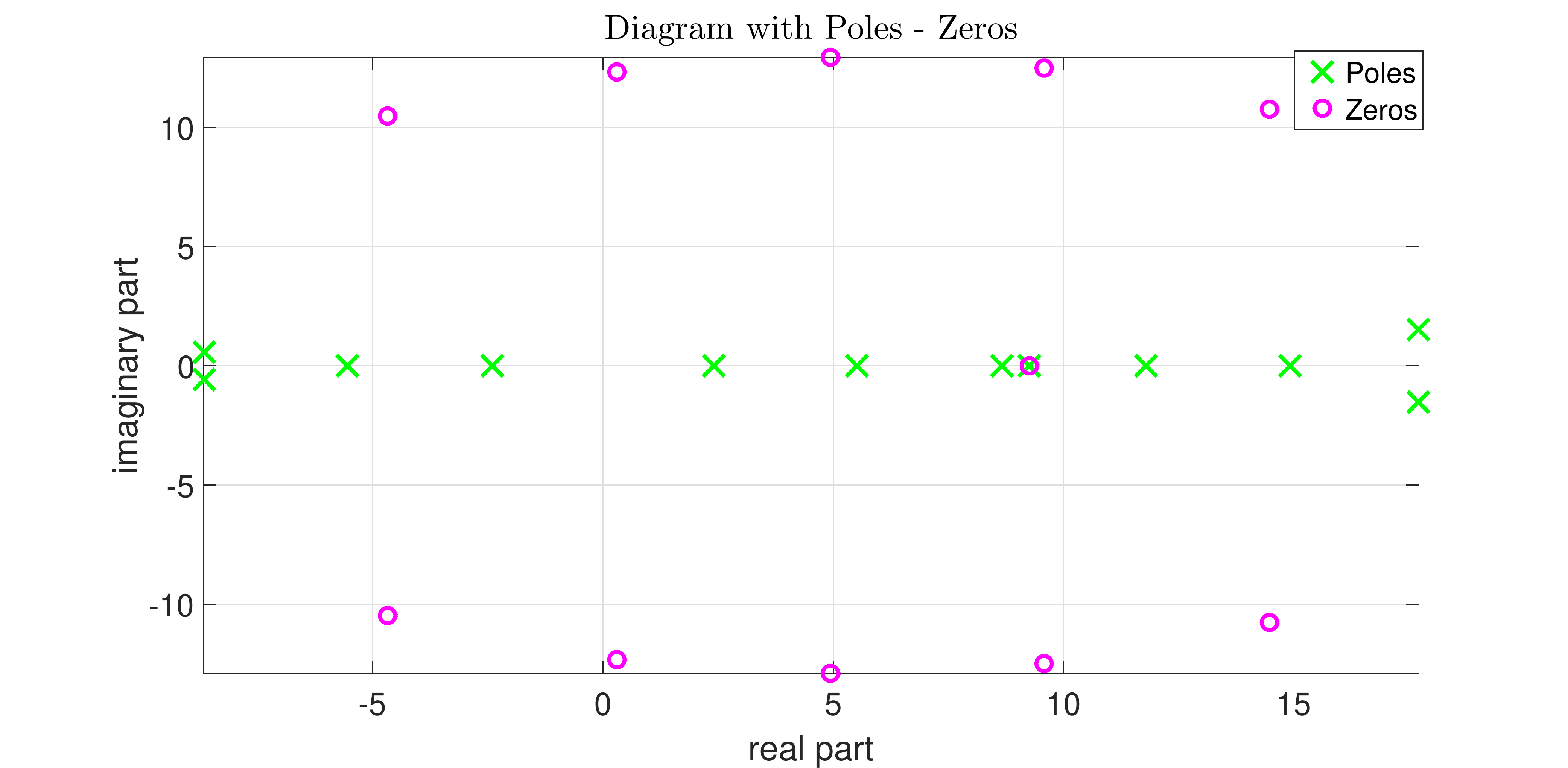}
  \caption{Poles and zeros diagram}
  \label{fig:sfig1}
\end{subfigure}%
\caption{Approximation results of the $H(s)$, over $\bOmega$, with the recursive Loewner algorithm approach.}
\label{fig:fig}
\end{figure}

In more detail the Figure 13a describes the interpolation points as shown in the right hand figure, section 2, page 5 with a depiction of these points on the corresponding samples of $H(s)$. The Figure 13b shows the 12th order approximant with the recursive Loewner framewoork approach superimposed on the plot of $H(s)$, over the dense grid $\bOmega_{grid}=[x_{1},...,x_{500}]\times[y_{1},...,y_{500}]\subset \bOmega$. As a result we see the error plot in Figure 13c which shows for each point in the $\bOmega_{grid}$, the absolute distance $|H_{r}(s)-H(s)|$, with $s \in \bOmega_{grid}$. The order of error is $O(10^{-10})$. Last, the corresponding poles/zeros are shown in the Figure 13d and for the exact values, someone can see Appendix 4 number 6.\\
\vspace{3mm}\\
\textbf{Remark 2.2.2.1} Since the 8th pole and the last zero above are equal up to $10^{-5}$, (Figure 13d or Appendix 4 number 6, they can not be eliminated because it is less than the error. So in this case recursive algorithm produce 12th order approximant with almost a pole/zero cancellation. 
\newpage
\subsubsection{The AAA Alorithm over Uniformly Distributed Points}

In this subsection we present the approximation with AAA algorithm where the input uses 2000 points as in figure 2b. After many realizations with uniformly distributed grid (2000 points), we noticed that the error varies between $O(1e-11)$, to $O(1e-13)$. We kept one of the best. 
\begin{figure}[h!]
\begin{subfigure}{.5\textwidth}
  \centering
  \includegraphics[width=1\linewidth]{sampling2000pointsBessel-eps-converted-to.pdf}
  \caption{Corresponding sampling points of $H(s)$.}
  \label{fig:sfig1}
\end{subfigure}%
\begin{subfigure}{.5\textwidth}
  \centering
  \includegraphics[width=1\linewidth]{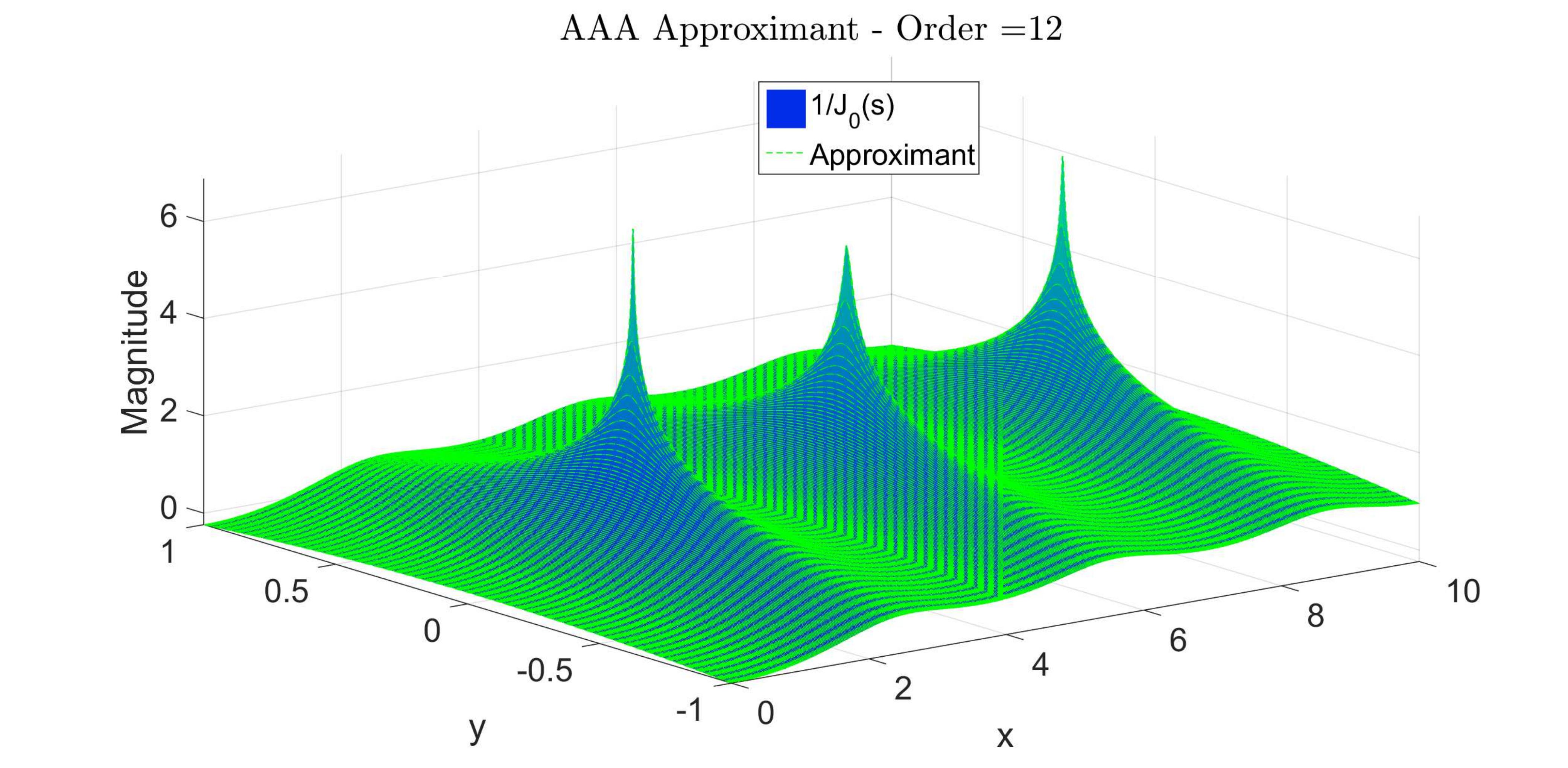}
  \caption{Approximant superimposed with $H(s)$, in $\bOmega$.}
  \label{fig:sfig2}
\end{subfigure}
\begin{subfigure}{.5\textwidth}
  \centering
  \includegraphics[width=1\linewidth]{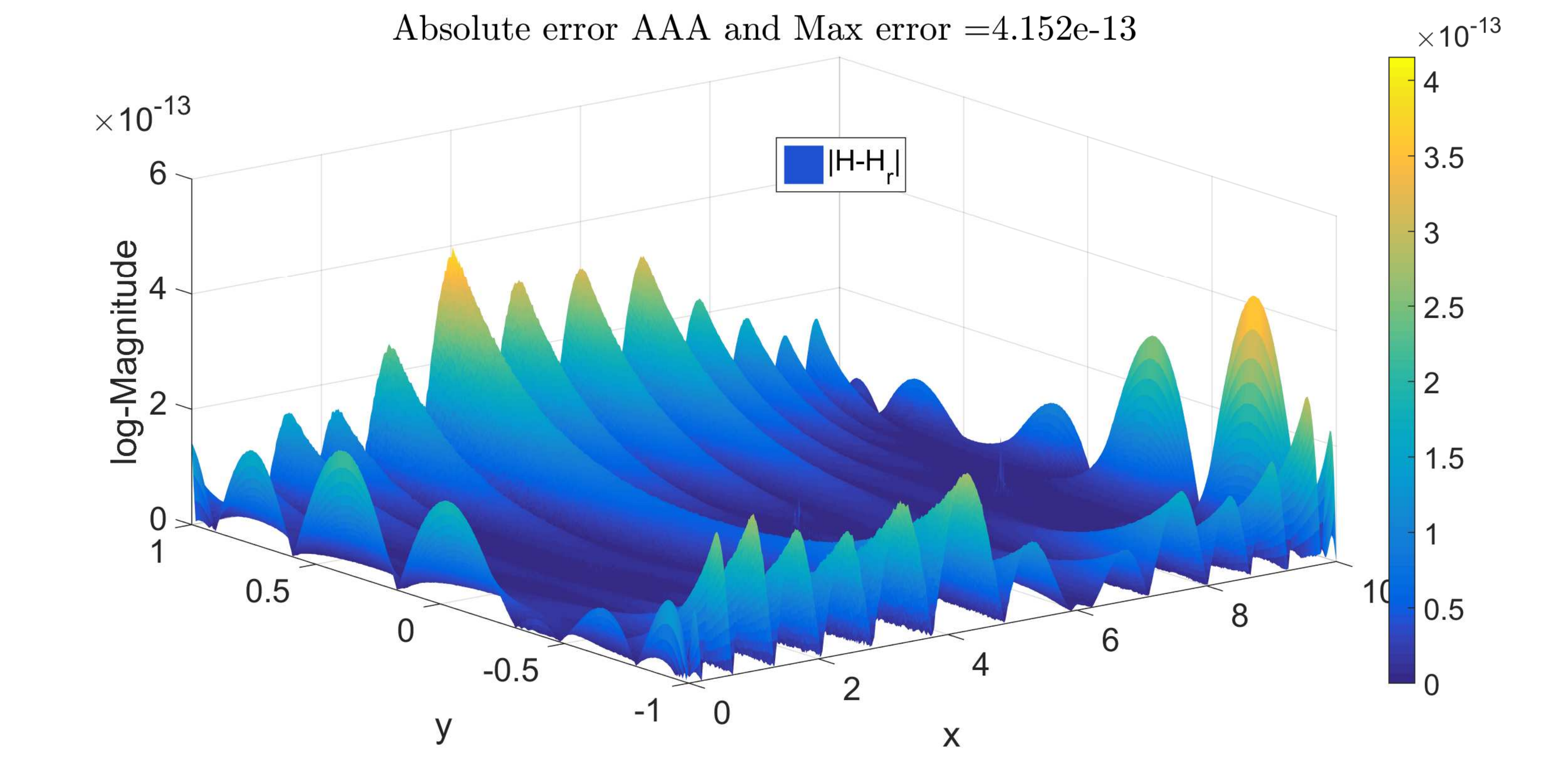}
  \caption{Error evaluation over the $\bOmega$.}
  \label{fig:sfig1}
\end{subfigure}%
\begin{subfigure}{.5\textwidth}
  \centering
  \includegraphics[width=1\linewidth]{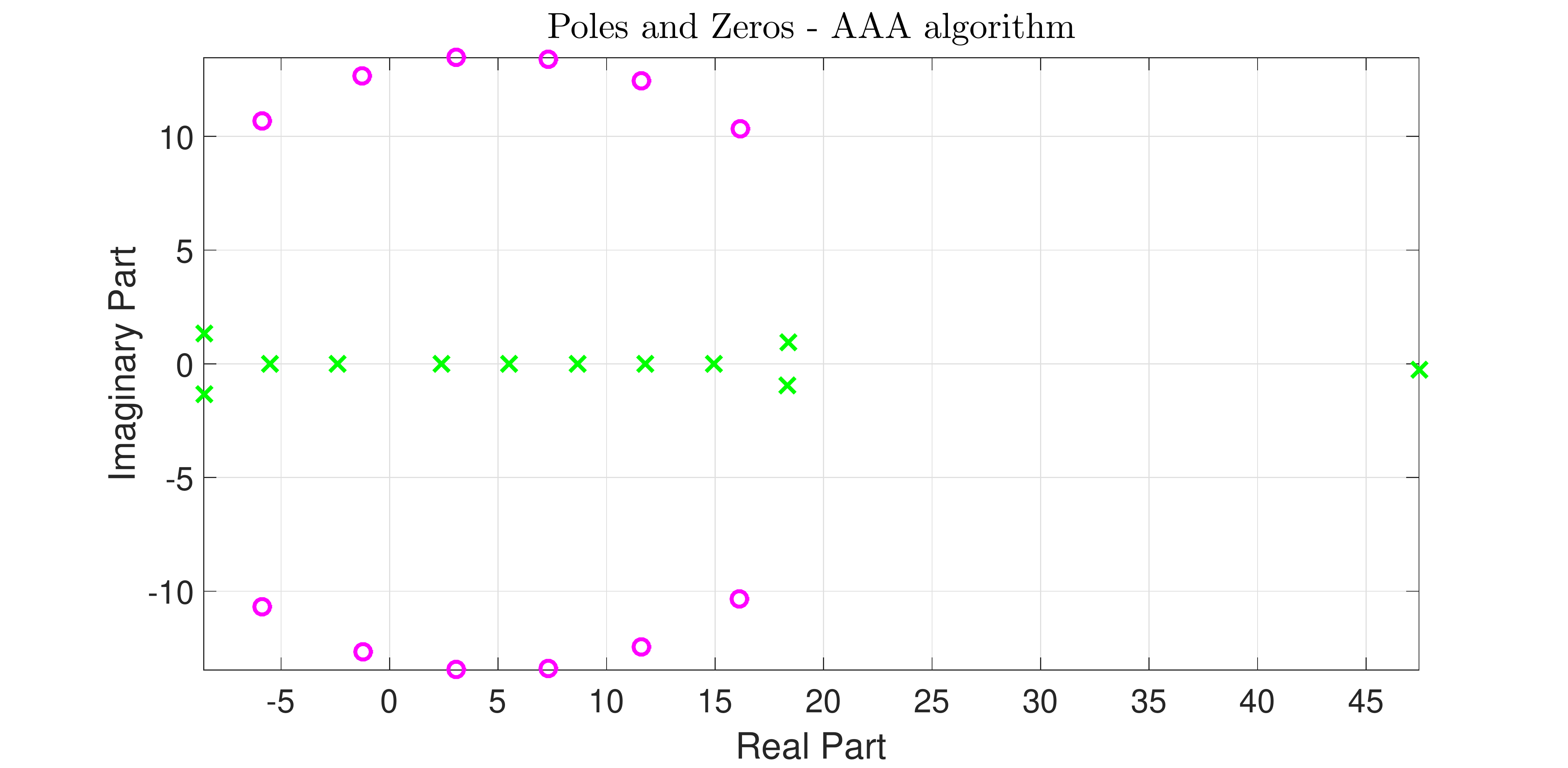}
  \caption{Poles and zeros diagram}
  \label{fig:sfig1}
\end{subfigure}%
\caption{Approximation results of the $H(s)$, over $\bOmega$, with the AAA algorithm approach.}
\label{fig:fig}
\end{figure}

In more detail the Figure 14a describes the interpolation points as shown in the right hand figure, section 2, page 5 with a depiction of these points on the corresponding samples of $H(s)$. The Figure 14b shows the 12th order approximant with the AAA algorithm approach superimposed on the plot of $H(s)$, over the dense grid $\bOmega_{grid}=[x_{1},...,x_{500}]\times[y_{1},...,y_{500}]\subset \bOmega$. As a result we see the error plot in Figure 14c which shows for each point in the $\bOmega_{grid}$, the absolute distance $|H_{r}(s)-H(s)|$, with $s \in \bOmega_{grid}$. The order of error is $O(10^{-13})$. Last, the corresponding poles/zeros are shown in the Figure 14d and for the exact values, someone can see Appendix 4 number 7.\\
\begin{figure}[h]
\centering
\includegraphics[width=0.7\linewidth, height=4cm]{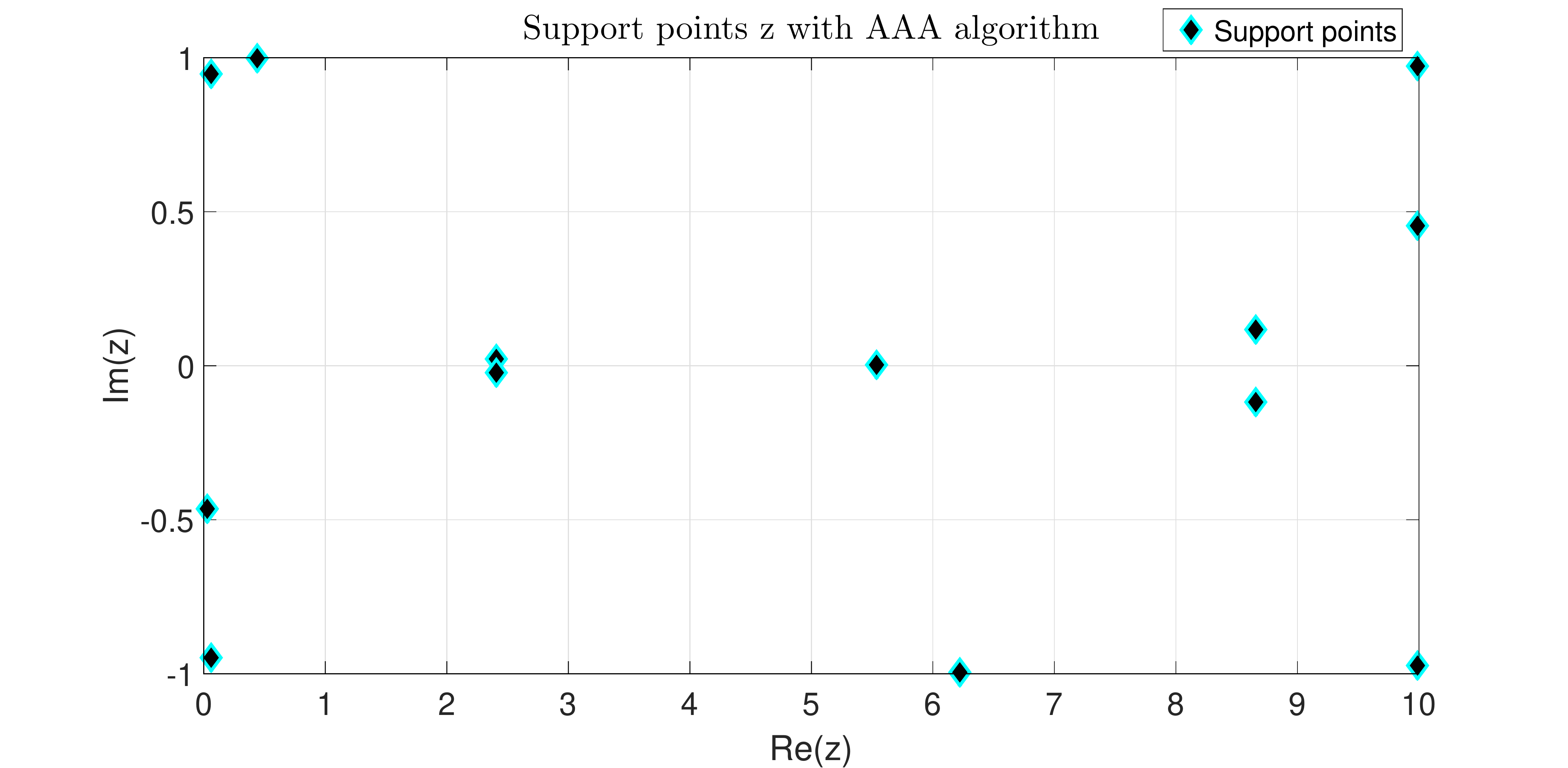} 
\caption{Support points from AAA algorithm in $\bOmega$.}
\end{figure}

\newpage
\subsubsection{Vector Fitting for uniformly distributed points}
In this subsection we present the approximation with VF algorithm where we use 2000 (conjugate pairs) uniformly distributed over the domain $\bOmega$, as in figure 2b.  
\begin{figure}[h!]
\begin{subfigure}{.5\textwidth}
  \centering
  \includegraphics[width=1\linewidth]{sampling2000pointsBessel-eps-converted-to.pdf}
  \caption{Corresponding sampling points of $H(s)$.}
  \label{fig:sfig1}
\end{subfigure}%
\begin{subfigure}{.5\textwidth}
  \centering
  \includegraphics[width=1\linewidth]{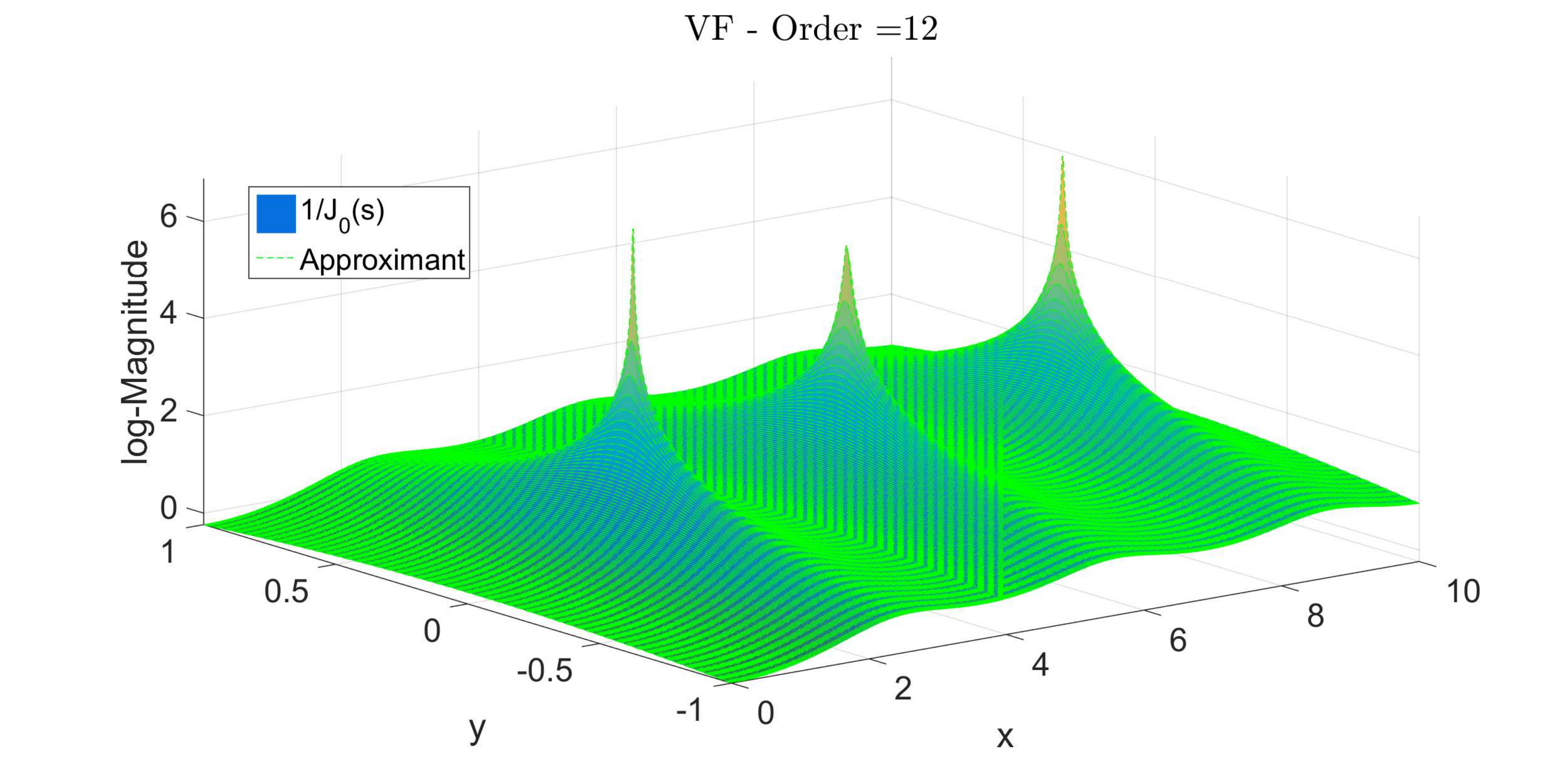}
  \caption{Approximant superimposed with $H(s)$, in $\bOmega$.}
  \label{fig:sfig2}
\end{subfigure}
\begin{subfigure}{.5\textwidth}
  \centering
  \includegraphics[width=1\linewidth]{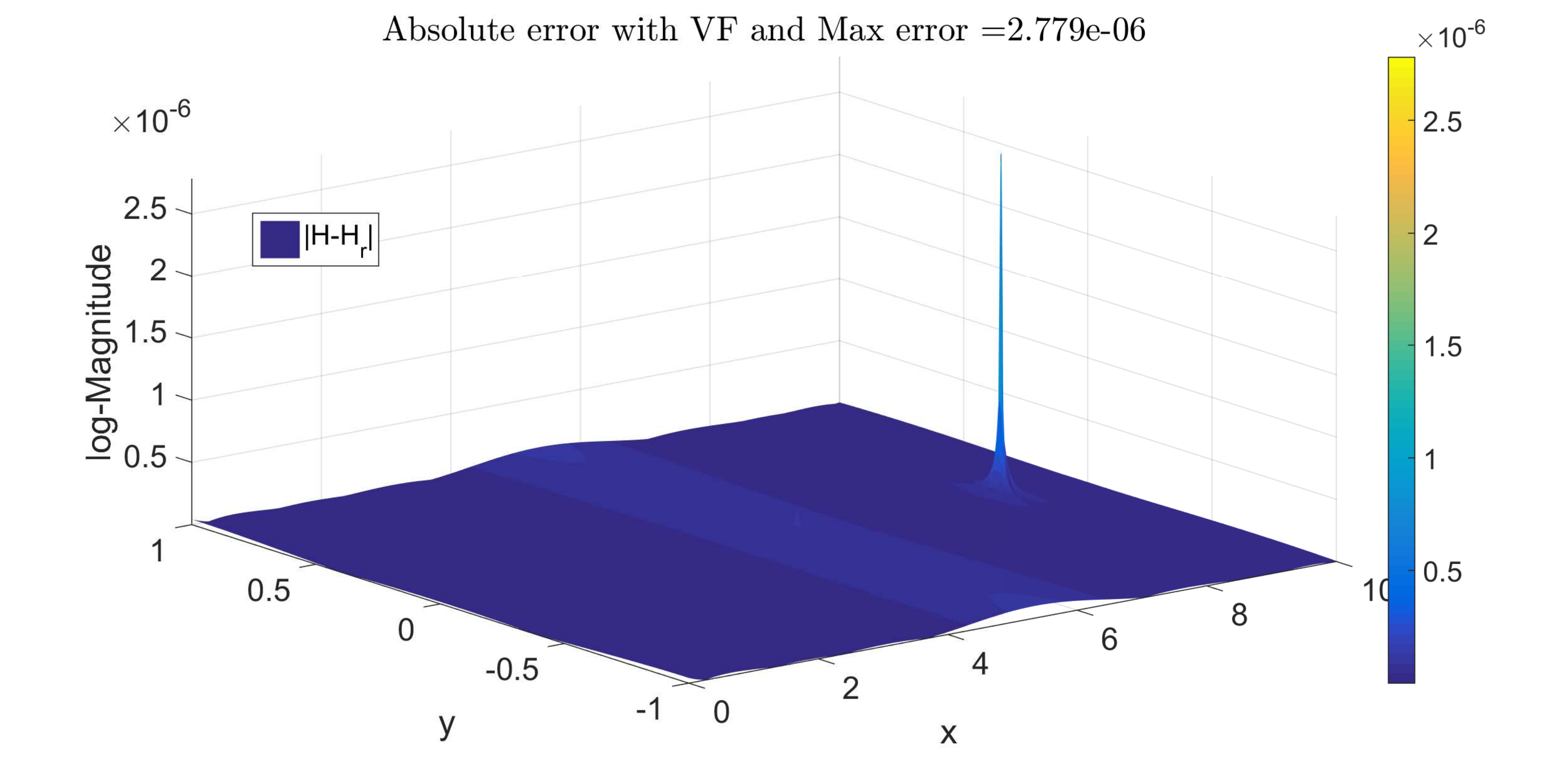}
  \caption{Error evaluation over the $\bOmega$.}
  \label{fig:sfig1}
\end{subfigure}%
\begin{subfigure}{.5\textwidth}
  \centering
  \includegraphics[width=1\linewidth]{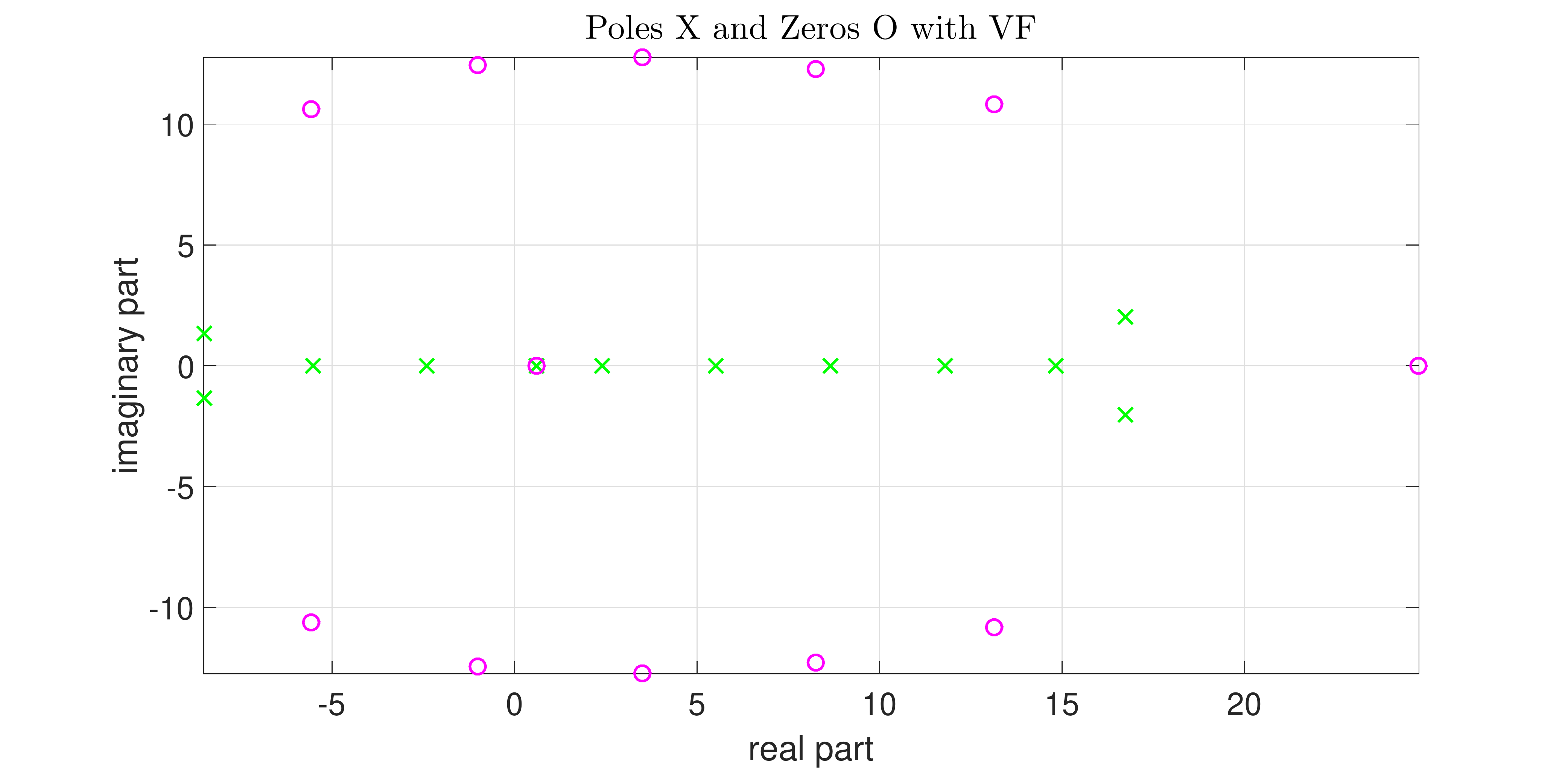}
  \caption{Poles and zeros diagram}
  \label{fig:sfig1}
\end{subfigure}%
\caption{Approximation results of the $H(s)$, over $\bOmega$, with the Vector Fitting algorithm approach.}
\label{fig:fig}
\end{figure}

In more detail the Figure 16a describes the interpolation points as shown in the right hand figure, section 2, page 5 with a depiction of these points on the corresponding samples of $H(s)$. The Figure 16b shows the 12th order approximant with the AAA algorithm approach superimposed on the plot of $H(s)$, over the dense grid $\bOmega_{grid}=[x_{1},...,x_{500}]\times[y_{1},...,y_{500}]\subset \bOmega$. As a result is the error plot in Figure 16c which shows for each point in the $\bOmega_{grid}$, the absolute distance $|H_{r}(s)-H(s)|$, with $s \in \bOmega_{grid}$. The order of error is $O(10^{-6})$. Last, the corresponding poles/zeros are shown in the Figure 16d and for the exact values, someone can see Appendix 4 number 8.\\
\vspace{2mm}\\
\textbf{Remark:} Since the one pole/zero pair above are equal up to $10^{-9}$, (Figure 16d) or Appendix 4 number 8, they can be eliminated. After this \textbf{zero/pole} cancellation we obtain the same order of accuracy with an approximate of order $r=11$.
\newpage
\section{Conclusion}
\subsection{Error comparison}
\begin{center}
  \begin{tabular}{ | c | c | c | c | c |}
    \hline
    Case/Method & Loewner & recursive Loewner & AAA & VF \\ \hline
    1st: 2121 structured points & $O(10^{-11})$ & $O(10^{-10})$ & $O(10^{-11})$ & $O(10^{-6})$ \\ \hline
    2nd: 2000 uniformly points & $O(10^{-11})$ & $O(10^{-10})$ & $O(10^{-13})$ & $O(10^{-6})$\\
    \hline
  \end{tabular}
\end{center}
\vspace{2mm}
The above results show that all methods constructed low order models with very good accuracy. Among them, the methods that built models with the highest accuracy are the direct method with the Loewner Framework and the iterative AAA algorithm. It is worthwhile to note that these two methods, in contrast to (recursive Loewner Framework and Vector Fitting), do not require the a-priori assignment of the order of the reduced system. The Loewner framework method resolves an entire SVD and from the way the singular values decay one can choose the truncation order. This gives a trade off between accuracy and complexity. On the other hand, the AAA algorithm, constructs the model by selecting on each step one more support point and solving a SVD on each iteration step. When the error drops down, AAA fullfills its exit criterion (i.e. tol=1e-13). Another point that is worth mentioning is that the Loewner framework method constructs models $(r-1,r)$,$r\in \mathbb{N}$, with real symmetry as opposed to the AAA algorithm which makes approximants of $(r,r)$, and with no real symmetry. The first number of the previous parenthesis stands for the order of the numerator and the second for the denominator order. This approach gives the opportunity for the Loewner framework method to introduce an arbitrary parameter D (parametrize the system, details see \cite{tutorial} \textbf{theorem 1.13 b}) where its calculation can reduce the error as in [7]. Another point is, that if we wanted to have the best possible approximation with the Loewner Framework, we could choose order 13th in both cases (Figure 2). The error would be: $O(10^{-13})$. We maintained a lower-order model (11th for Figure 2a and 12th for Figure 2b) in order to do a comparison with the AAA algorithm. 

Last but not least, as we feed the system with conjugate measurements in both cases, someone should be able to reveal the underlying real model. The AAA algorithm does not impose real symmetry. For all the above methods we build a real approximant except for the AAA aproximant. Introducing some changes to the AAA alorithm, now we are able to construct a reduced model with real symmetry. In this particular case, we noticed that the order of the approximant becomes higher and produce many poles/zeros cancellations (Froissart doublets) up to precision. In order to tackle this problem, it is needed to enable the function "cleanup" as explained in \cite{AAA}. Both functions increase the complexity for the construction of a real approximant.

\subsection{Complexity comparison}
\textbf{The Complexity of the Loewner Framework algorith and AAA algorithm.}\\
The main computational cost, for both methods comes from SVD calculation. For a matrix A with dimension $m\times n$, the SVD cost is $O(m^{2}n+mn^{2}+n^3)$. Since the Loewner Framework involves SVDs of dimension $(\frac{N}{2}\times \frac{N}{2})$, the resulting complexity is $2O(\frac{3N^3}{8})\approx O(N^3)$, flops.
Algorithm AAA involves SVDs of dimensions $(M-j)\times j$, with $j=1,2,...,m$, ($m$, stands for the order) its complexity is $O(Mm^3)$, flops \cite{AAA} instead of the Loewner framework which solves only one SVD with complexity $O(M^3)$. For this problem on Bessel approximation, two methods succeed to construct reduced models with low order approximants. The Loewner framework would produce better computational results in the case that we have a slow decay of the singular values, so it would be needed to truncate in a higher order. This is something that will increase the computational cost of AAA algorithm. 


\section{Appendix}
\textbf{The Bessel function zeros of the first kind}.
The first 6 solutions of the equation $J_{0}(s)=0$, with 15 digits precision are:
\begin{center}
\textbf{Bessel original Zeros with 15 digits}
\tiny
$=\left(\begin{array}{c}
\textbf{2.40482555769577}\\ \textbf{5.52007811028631}\\ \textbf{8.65372791291101}\\11.7915344390142 \\14.9309177084877\\18.0710639679109\end{array}\right)$.
\end{center}
\begin{enumerate}
\item Poles/zeros from the Loewner Framework with 2121 points grid
\begin{center}
\tiny
$\textbf{Poles}=\left(\begin{array}{c} -8.32213293322054- 1.4252\, \mathrm{i}\\ -8.32213289862456+ 1.4252\, \mathrm{i}\\ -5.51461491999547\\ -2.40481847965605\\ \textbf{2.40482555769577}\\ \textbf{5.52007811028631}\\ \textbf{8.65372791291101}\\ 11.7915356008908\\ 14.9135964357538\\ 17.6548692348549- 1.561\, \mathrm{i}\\ 17.654869354827+ 1.561\, \mathrm{i} \end{array}\right)$,
$\textbf{Zeros}=\left(\begin{array}{c} -4.8491 - 10.766\, \mathrm{i}\\ -4.8491 + 10.766\, \mathrm{i}\\ 0.12013 - 12.537\, \mathrm{i}\\ 0.12013 + 12.537\, \mathrm{i}\\ 4.785 + 13.066\, \mathrm{i}\\ 4.785 - 13.066\, \mathrm{i}\\ 9.4384 + 12.591\, \mathrm{i}\\ 9.4384 - 12.591\, \mathrm{i}\\ 14.367 + 10.868\, \mathrm{i}\\ 14.367 - 10.868\, \mathrm{i} \end{array}\right)$.
\end{center}

\textbf{Remark:} With \textbf{bold} are the poles inside the $\bOmega$, domain.
\item Poles/zeros from the recursive Loewner approach with 2121 points grid
\begin{center}
\tiny
$PolesRec=\left(\begin{array}{c} -8.64763053774167- 0.874672\, \mathrm{i}\\ -8.64763053774167+ 0.874672\, \mathrm{i}\\ -5.53736536596127\\ -2.40481506507015\\ \textbf{2.40482555769639}\\\textbf{5.52007811028556}\\ \textbf{8.65372791291091}\\ 11.791534687103\\ 14.9207741074572\\ 17.7972322878303- 1.42773\, \mathrm{i}\\ 17.7972322878303+ 1.42773\, \mathrm{i} \end{array}\right)$.
$ZerosRec=\left(\begin{array}{c} -4.81774 - 10.3428\, \mathrm{i}\\ -4.81774 + 10.3428\, \mathrm{i}\\ 0.196012 - 12.134\, \mathrm{i}\\ 0.196012 + 12.134\, \mathrm{i}\\ 4.89894 - 12.6902\, \mathrm{i}\\ 4.89894 + 12.6902\, \mathrm{i}\\ 9.58028 - 12.2575\, \mathrm{i}\\ 9.58028 + 12.2575\, \mathrm{i}\\ 14.5169 - 10.5921\, \mathrm{i}\\ 14.5169 + 10.5921\, \mathrm{i} \end{array}\right)$.\\
\end{center}
\item Poles/zeros from the AAA approach with 2121 points grid
\begin{center}
\tiny
$\text{PolesAAA}=\left(\begin{array}{c} -7.70989384927742+ 1.98333\, \mathrm{i}\\ -7.64706106530066- 1.84221\, \mathrm{i}\\ -5.46205580369784 + 0.00647\,\mathrm{i}\\ -2.40483752963882\\ \textbf{2.40482555769577}\\ \textbf{5.52007811028632}\\ \textbf{8.65372791291101}\\ 11.7915347392982\\ 14.9221430048388+0.00059\,\mathrm{i}\\ 17.8709511796541- 1.35084\, \mathrm{i}\\ 17.9082708199292+1.39389\, \mathrm{i} \end{array}\right)$,
$\text{ZerosAAA}=\left(\begin{array}{c} -4.05365 - 10.9078\, \mathrm{i}\\ -3.99775 + 11.0922\, \mathrm{i}\\ 0.832526 - 12.5068\, \mathrm{i}\\ 0.905811 + 12.6642\, \mathrm{i}\\ 5.42377 - 12.9188\, \mathrm{i}\\ 5.49875 + 13.0562\, \mathrm{i}\\ 9.99811 - 12.3601\, \mathrm{i}\\ 10.0689 + 12.476\, \mathrm{i}\\ 14.8259 - 10.5757\, \mathrm{i}\\ 14.892 + 10.666\, \mathrm{i}\\ -22.5906 + 2.52206\, \mathrm{i} \end{array}\right)$.\\
\end{center}
\textbf{Remark: Poles and Zeros do not appear as conjugates.}
\item Poles/zeros from the Vector Fitting algorithm approach with 2121 points grid
\begin{center}
\tiny
$\text{PolesVF}=\left(\begin{array}{c} -8.49629987353267- 1.3597\, \mathrm{i}\\ -8.49629987353267+ 1.3597\, \mathrm{i}\\ -5.51126493255679\\ -2.40482600718294\\ \color{red}{-0.157905556703521}\\ \textbf{2.40482555759485}\\ \textbf{5.52007811166818}\\ \textbf{8.65372790279672}\\ 11.7915516528236\\ 14.8446346140717\\ 16.8161497628161- 2.0178\, \mathrm{i}\\ 16.8161497628161+ 2.0178\, \mathrm{i} \end{array}\right)$,
$\text{ZerosVF}=\left(\begin{array}{c} \color{red}{-0.157905555129157}\, \\ -5.6578 - 10.624\, \mathrm{i}\\ -5.6578 + 10.624\, \mathrm{i}\\ -0.9698 + 12.28\, \mathrm{i}\\ -0.9698 - 12.28\, \mathrm{i}\\ 3.683 - 12.71\, \mathrm{i}\\ 3.683 + 12.71\, \mathrm{i}\\ 8.3717 + 12.366\, \mathrm{i}\\ 8.3717 - 12.366\, \mathrm{i}\\ 13.196 - 10.889\, \mathrm{i}\\ 13.196 + 10.889\, \mathrm{i}\\ 25.816\end{array}\right).$
\end{center}

\textbf{Remark: red color pole/zero pair shows cancellation}
\item Poles/zeros from the Loewner Framework with 2000 uniformly points
\begin{center}
\tiny
$\text{PolesRL}=\left(\begin{array}{c} -8.05158320855723- 1.7303\, \mathrm{i}\\ -8.05158320855723+ 1.7303\, \mathrm{i}\\ -5.48778972998639\\ -2.40483192639586\\ \textbf{2.40482555769577}\\\textbf{5.52007811028632}\\ \textbf{8.65372791291101}\\ 11.7915339297017\\ 15.0015209902362 - 0.222\, \mathrm{i}\\ 15.0015209902362 + 0.222\, \mathrm{i}\\ 17.377289961194- 2.03\, \mathrm{i}\\ 17.377289961194+ 2.03\, \mathrm{i} \end{array}\right)$,
$\text{ZerosRL}=\left(\begin{array}{c} -4.7536 - 11.207\, \mathrm{i}\\ -4.7536 + 11.207\, \mathrm{i}\\ 0.17552 - 12.998\, \mathrm{i}\\ 0.17552 + 12.998\, \mathrm{i}\\ 4.8076 - 13.534\, \mathrm{i}\\ 4.8076 + 13.534\, \mathrm{i}\\ 15.242\\ 9.4292 - 13.06\, \mathrm{i}\\ 9.4292 + 13.06\, \mathrm{i}\\ 14.319 - 11.333\, \mathrm{i}\\ 14.319 + 11.333\, \mathrm{i}\\  \end{array}\right)$.
\end{center}
\item Poles/zeros from the recursive Loewner with 2000 uniformly points
\begin{center}
\tiny
$polesREC=\left(\begin{array}{c} -8.66473644838372- 0.55511\, \mathrm{i}\\ -8.66473644838372+ 0.55511\, \mathrm{i}\\ -5.55846358891533\\ -2.40478725513252\\ 2.40482555769579\\ 5.52007811028631\\ 8.65372791291094\\ \color{red}{9.2576}\color{black}3698485788\\ 11.7915351867464\\ 14.9164167045242\\ 17.7116227515344- 1.513\, \mathrm{i}\\ 17.7116227515344+ 1.513\, \mathrm{i} \end{array}\right)$,
$ZerosREC=\left(\begin{array}{c}\ 14.465 + 10.772\, \mathrm{i}\\ 14.465 - 10.772\, \mathrm{i}\\ -4.6782 + 10.478\, \mathrm{i}\\ -4.6782 - 10.478\, \mathrm{i}\\ 9.5668 + 12.474\, \mathrm{i}\\ 9.5668 - 12.474\, \mathrm{i}\\ 0.29543 + 12.334\, \mathrm{i}\\ 0.29543 - 12.334\, \mathrm{i}\\ 4.9411 + 12.915\, \mathrm{i}\\ 4.9411 - 12.915\, \mathrm{i}\\ \color{red}{9.2576}\color{black}4464786048\ \end{array}\right)$.
\end{center}
\item Poles/zeros from the AAA with 2000 uniformly points
\begin{center}
\tiny
$\text{PolesAAA}=\left(\begin{array}{c} -8.56002747182546- 1.3357\, \mathrm{i}\\ -8.54868176358933+ 1.3357\, \mathrm{i}\\ -5.51154958824584\\ -2.40482617183982\\ \textbf{2.40482555769577}\\ \textbf{5.52007811028630}\\ \textbf{8.65372791291101}\\ 11.7915344794608\\ 14.9286212759068+ 0.000044838\, \mathrm{i}\\ 18.349207798415518.349 - 0.93836\, \mathrm{i}\\ 18.3564686110828+ 0.94294\, \mathrm{i}\\ 47.4348759491114- 0.25113\, \mathrm{i}  \end{array}\right)$,
$\text{ZerosAAA}=\left(\begin{array}{c} -5.8598 - 10.679\, \mathrm{i}\\ -5.8806 + 10.684\, \mathrm{i}\\ -1.2407 - 12.66\, \mathrm{i}\\ -1.2625 + 12.659\, \mathrm{i}\\ 3.062 + 13.453\, \mathrm{i}\\ 3.0802 - 13.46\, \mathrm{i}\\ 7.3122 + 13.387\, \mathrm{i}\\ 7.3233 - 13.396\, \mathrm{i}\\ 11.605 + 12.456\, \mathrm{i}\\ 11.608 - 12.463\, \mathrm{i}\\ 16.139 - 10.328\, \mathrm{i}\\ 16.142 + 10.327\, \mathrm{i} \end{array}\right)$.
\end{center}
\textbf{Remark: Poles and Zeros do not appear as conjugates.}
\item Poles/zeros from the VF with 2000 uniformly points
\begin{center}
\tiny
$\text{PolesVF}=\left(\begin{array}{c} -8.51428135857578- 1.3424\, \mathrm{i}\\ -8.51428135857578+ 1.3424\, \mathrm{i}\\ -5.51184166025835\\ -2.40482593907193\\ \color{red}{0.598123905184851}\\ \textbf{2.40482555782458}\\ \textbf{5.52007810888635}\\ \textbf{8.65372792599631}\\ 11.7915548184731\\ 14.8352178424186\\ 16.7433767283711- 2.0413\, \mathrm{i}\\ 16.7433767283711+ 2.0413\, \mathrm{i} \end{array}\right)$,
$\text{ZerosVF}=\left(\begin{array}{c} \color{red}{0.598123905506776} \\\ -5.5835 - 10.635\, \mathrm{i}\\ -5.5835 + 10.635\, \mathrm{i}\\ -1.0088 + 12.445\, \mathrm{i}\\ -1.0088 - 12.445\, \mathrm{i}\\ 3.4938 + 12.744\, \mathrm{i}\\ 3.4938 - 12.744\, \mathrm{i}\\ 8.2407 - 12.29\, \mathrm{i}\\ 8.2407 + 12.29\, \mathrm{i}\\ 13.127 + 10.826\, \mathrm{i}\\ 13.127 - 10.826\, \mathrm{i}\\ 24.775 \, \ \end{array}\right)$.
\end{center}
\textbf{Remark: red color pole/zero pair shows cancellation}
\end{enumerate}
\newpage
\footnotesize
\medskip
\bibliographystyle{siam}

\end{document}